%% file: DoesElasticStressModifyEquilbriumCornerAngle.tex
\documentclass[3p,letter]{elsarticle}
\usepackage{lineno,hyperref}
\usepackage{amsmath,bm}
\modulolinenumbers[5]
\journal{J Mech Phys Solids}

\usepackage[ruled,vlined]{algorithm2e}
\usepackage{subcaption}
\newcommand{\hs}{\hspace{5mm}}

\begin{document}

\begin{frontmatter}
\title{Does elastic stress modify the equilibrium corner angle? \tnoteref{mytitlenote}}


\author[mymainaddress,presentaddress]{Weiqi Wang}
\cortext[mycorrespondingauthor]{Corresponding author}
\ead{weiqi.wang@concordia.ca}

\author[mymainaddress]{Brian J. Spencer\corref{mycorrespondingauthor}}
\ead{spencerb@buffalo.edu}

\address[mymainaddress]{Department of Mathematics, University at Buffalo, Buffalo, NY 14260, USA}
\address[presentaddress]{Department of Mathematics and Statistics,
Concordia University, Montreal, Quebec, Canada H3G1M8}

\begin{abstract}
\input{abstract.tex}
\end{abstract}

\begin{keyword} \small
corner stress singularity \sep
voids and inclusions \sep
variational calculus \sep
boundary integral equations \sep
asymptotic analysis 
\end{keyword}

\end{frontmatter}


\input{Sec1_Introduction}

\input{Sec2_Mathematical_Formulation}

\input{Sec3_Numerical_Method}

\input{Sec4_Result}

\input{Sec5_Discussion}

\input{Sec6_Conclusion}

\appendix

\input{Appendix_A}

\input{Appendix_B}

\input{Appendix_C}

\input{Appendix_D}

\section*{Funding} 
This work was supported by a grant from the Simons Foundation (Award \#354717, BJS).

\bibliography{mybibfile}

\end{document}

%% file: abstract.tex
We consider the influence of elasticity and anisotropic surface energy on the energy-minimizing shape of a two-dimensional void under biaxial loading.
In particular, we consider void shapes with corners for which the strain energy density is singular at the corner.  
The elasticity problem is formulated as a boundary integral equation using complex potentials.  By incorporating the asymptotic behavior of the singular elastic fields at corners of the void, we develop a numerical spectral method for determining the stress for a class of arbitrary void shapes and corner angles.  
We minimize the total energy of surface energy and elastic potential energy using calculus of variations to obtain an Euler-Lagrange equation on the boundary that is coupled to the elastic field.  The shape of the void boundary is determined using a numerical spectral method that simultaneously determines the equilibrium void shape and singular elastic fields. 
Our results show that the precise corner angles that minimize the total energy are not affected by the presence of the singular elastic fields. However, the stress singularity on the void surface at the corner must be balanced by a singularity in the curvature at the corner that effectively changes the macroscopic geometry of the shape and effectively changes the apparent corner angle.   
These results reconcile the apparent contradiction regarding the effect
of elasticity on equilibrium corner angles in the results of Srolovitz and Davis (2001) and Siegel, Miksis and Voorhees (2004), and identify an important nontrivial singular behavior associated with corners on 
free-boundary elasticity problems.

%% file: Sec1_Introduction.tex
\section{Introduction}\label{Sec:1}

Microstructures in materials are typically driven to energy-minimizing configurations.  For multi-phase systems the microstructure is characterized by the boundary between the phases. Mathematically, the microstructure can be viewed as a free boundary problem for the shape or geometry of the phases that corresponds to minimizing the energy of the system with respect to volumetric and surface energy contributions.

In this paper we are interested in microstructure shapes that involve corner angles or contact angles between phases.  
 There are many factors that can affect the contact angle, 
such as surface properties of the solid, surface roughness \cite{wenzel1936}, and the temperature and pressure of the system \cite{young1805}. 
In this work we investigate the effect of elasticity on corner angles.
In particular, we consider a prototype problem of a void in a two-dimensional (2D) solid 
to isolate the effect of 
anisotropic surface energy  \cite{choi2009} and elastic stress \cite{srolovitz1989} on
 the corner angle.

Because of the surface energy anisotropy, the equilibrium void shape may exclude surface orientations with high surface energy, giving rise to corners on the equilibrium shape \cite{cabrera1964,Voorhees1984}. 
The corner angles for the equilibrium shapes in the absence of elastic stresses are given by the Wulff shape \cite{wulff1901} and are well known (see for example \cite{spencer2004} and references therein). 
With the presence of the elastic stress, a corner angle can generate an elastic singularity.  For example, in the wedge geometry (i.e., a wedge-shaped elastic solid in a plane) with a solid corner angle greater than $\pi$ (void angle less than $\pi$) may have an elastic stress singularity at the vertex of the wedge \cite{williams1952stress}.
 Therefore, for the anisotropic void shape 
with corners, the stress field may have corner singularities due to the loading of the solid.  Thus the elastic energy may have a singularity at the void corner which could affect the total energy 
(surface energy and elastic energy) for a given void shape. 

Since the equilibrium corner angles are determined by minimization of the total energy with respect to the entire void shape, a natural question is whether the presence of the stress singularity will modify the corner angles.
Srolovitz and Davis \cite{srolovitz2001} analyze the elastic energy locally near the corner and use scaling arguments to conclude that stresses do
not modify the corner angle from the no-stress results. In contrast, Siegel, Miksis and Voorhees \cite{siegel2004} 
consider the energy-minimization problem for the shape of a void in a stressed solid by using a corner energy regularization term in which there
is a curvature term as an energy penalty for corners.  By considering the case of small regularization, the numerical results appear to converge to sharp-corner shapes. However,
the apparent corner angles show a distinct dependence on the amount of elastic loading (see Fig.\,7 in \cite{siegel2004}) which suggests that the corner angle does depend on elastic stress.  Thus, there is an apparent contradiction between Srolovitz and Davis \cite{srolovitz2001} and Siegel et al \cite{siegel2004} about whether elasticity causes a change to the corner angle of a void or not.

Motivated by the apparent contradiction on the effect of elasticity on  the corner angle of a void, 
in this paper we develop a numerical method to find the equilibrium void shape with corners and elastic singularities.  As is known, the corner singularity can cause difficulties and instability in
numerical methods especially in the finite element method \cite{oh1995method,liu2012singular}.  The work of this paper develops a boundary integral formulation with a piecewise-continuous spectral method to solve the elasticity
problem for a general class of void shapes with corner singularity.   The behavior of the elastic singularity is incorporated directly into the numerical method using asymptotic analysis of the corner region.  

The elasticity solver is then used as a component in determining the shape of the void (including corner angles) that minimize the total energy.
Away from the corner, calculus of variations gives the classical Euler-Lagrange equation for the shape of the void which involves a balance of surface energy and elastic energy contributions.  Because of the elastic singularity at the corner, equilibrium also requires a singularity in the curvature which we determine from asymptotic analysis in the corner region and incorporate into our spectral representation for the void shape.

Our numerical method thus determines the equilibrium shape, corner angles and elastic field associated with a void. Our results confirm that the precise corner angles are not affected by elastic stress, in agreement with Srolovitz and Davis \cite{srolovitz2001}. However, the curvature singularity, which is driven by the elastic singularity, results in a change in the apparent corner angle from its stress-free value when the void shape is viewed macroscopically.  These results agree with the results shown in \cite{siegel2004} for the regularized-corner calculations.   Thus, our work  resolves the apparent contradiction between the results in \cite{srolovitz2001} and \cite{siegel2004} about the role of stress on corner angles but it also makes clear that the behavior in the corner region is singular and nontrivial.

The organization of the paper is as follows. Sec.\ 2 presents the mathematical formulation for the free boundary elasticity problem corresponding to energy minimization.  In Sec.\ 3 we describe the numerical method and its validation.  In Sec.\ 4 we obtain the results for the void shape and corner angle and highlight the change in the apparent corner angle due to elastic effects.  In Sec.\ 5 we discuss the parameters that control the change in the apparent corner angle from the perspective of void size and elastic loading.  Sec.\ 6 contains a concluding summary.

%% file: Sec2_Mathematical_Formulation.tex
\section{Mathematical Formulation} \label{Sec:2}

\begin{figure}[!t]
\centering
\includegraphics[width=75mm]{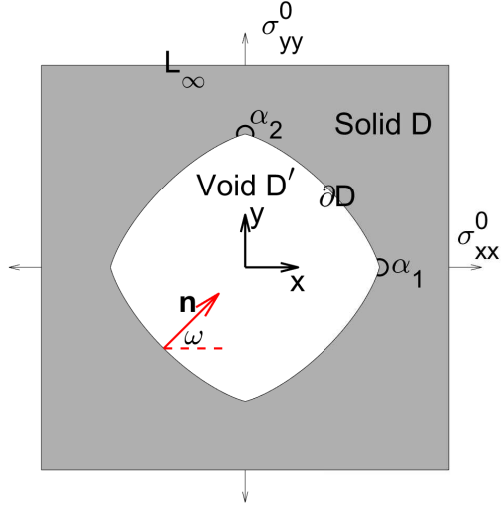}
\caption{Void geometry and notation.} \label{fig:1}
\end{figure}

In this section we describe the elasticity problem for the void and the equations for the void shape that correspond to minimizing the potential energy at fixed void area.   We consider the two-dimensional void geometry in the $(x,y)$ plane shown in Fig.~\ref{fig:1} where the origin is taken as the geometric center of the void, $D$ is the solid region exterior to the void, $D'$ is the void region, $\partial D$ is the piecewise-smooth void boundary with exterior unit normal $\bm{n}$ and orientation angle $\omega$, $\alpha_1, \alpha_2$ are solid corner angles on the $x,y$ axis, and 
$\sigma^0_{xx}$, $\sigma^0_{yy}$ are the applied biaxial stresses at infinity.

\subsection{Elasticity and elastic energy}
We denote the  the stress tensor as
\begin{equation}
\bm{\sigma} =\left[\begin{array}{cc}
\sigma_{xx} & \sigma_{xy}\\
\sigma_{xy} & \sigma_{yy}
\end{array}\right].
\end{equation}
We use  indicial subscript notation for the coordinates $(x,y) = (x_1,x_2)$ where repeated indices imply summation (e.g. $\sigma_{kk} = \sigma_{11}+\sigma_{22}$), and where subscripted commas denote differentiation with respect to the variables implied by indices after the comma (e.g. $\sigma_{12,2} = \partial \sigma_{12} /\partial x_2$).
The equations for mechanical equilibrium are
\begin{eqnarray}
\sigma_{ij,j} = 0 &&\hs \mbox{in $D$}, \label{eq:elasticity}\\
\sigma_{ij} n_j = 0 &&\hs \mbox{on $\partial D$}, \label{eq:3}\\
\sigma_{ij} \rightarrow \sigma^{\infty}_{ij} &&\hs \mbox{as $\sqrt{x^2 + y^2} \rightarrow \infty$},
\end{eqnarray}
where
\begin{equation}
\bm{\sigma^\infty} =\left[\begin{array}{cc}
\sigma_{xx}^0 & 0\\
0 & \sigma_{yy}^0
\end{array}\right]
\end{equation}
prescribes the biaxial stress applied at infinity.  
Defining the displacement vector $\bm{u} = (u_x, u_y)$ and the linearized strain tensor with components
\begin{equation}
\epsilon_{ij} = \tfrac{1}{2} ( u_{i,j} + u_{j,i} ),
\end{equation}
the compatibility condition for strain is 
\begin{equation}
\epsilon_{11,22} - 2 \epsilon_{12,12} + \epsilon_{22,11}=0 
\end{equation}
and
the constitutive law for isotropic plane-strain linear elasticity is
\begin{equation}
\sigma_{ij} = \frac{E}{1+\nu} \left[ \epsilon_{ij} + \frac{\nu}{1-2\nu} \delta_{ij}\epsilon_{kk} \right] \label{eq:constitutive_law}
\end{equation}
where $E$ is Young's modulus, $\nu$ is Poisson's ratio, and $\delta_{ij}$ is Kronecker delta.

Later we will derive the equations for the void shape that satisfy the free boundary conditions for minimization of total potential energy.  
For the elastic component of the potential energy, the applied biaxial stress means that the system is not isolated, so the elastic potential energy $\Pi_{\mbox{el}}$ is given by the total internal energy $E_{\mbox{int}}$ minus the work done by external forces $W_{\mbox{ext}}$ \cite{soutas2012potential},
\begin{equation}
\Pi_{\mbox{el}} = E_{\mbox{int}} - W_{\mbox{ext}},
\label{eq:elastic_potential_energy}
\end{equation}
where
\begin{equation}
E_{\mbox{int}} = \int \!\int_D \tfrac{1}{2}  \sigma_{ij} \epsilon_{ij} \, dA
\label{eq:elastic_internal_energy}
\end{equation}
and
\begin{equation}
W_{\mbox{ext}} = \int_{L_\infty} n_i \sigma_{ij} u_j \, ds
\label{eq:elastic_work}
\end{equation}
where $L_\infty$ is the circular contour described by $R = \sqrt{x^2+y^2}$ as $R \rightarrow \infty$ and $s$ is arclength.  In \ref{Appendix:A} we show following \cite{suo1994} that $E_{\mbox{int}} = \tfrac{1}{2} W_{\mbox{ext}}$ for the void and that 
\begin{equation}
\Pi_{\mbox{el}} = - \int \!\int_D \tfrac{1}{2} \sigma_{ij} \epsilon_{ij} \, dA= \tfrac{1}{2} \int_{\partial D} n_i \sigma_{ij}^\infty u_j \, ds - C \label{eq:12}
\end{equation}
where $C$ is a constant independent of the void shape.

\subsection{Void shape and surface energy}
We consider a void shape described in polar coordinates by $R \leq r(\theta)$, where $\theta$ is the polar angle measured counterclockwise from the positive $x$ axis and $R$ is the radial coordinate.  
The shape of the void is a free boundary problem for which the equilibrium shape minimizes the total potential energy, comprised of surface energy and elastic potential energy.   The surface energy is 
\begin{equation}
E_{\mbox{surf}} = \int_{\partial D} \gamma(\omega) \, ds
\end{equation}
where $\gamma(\omega)$ is the anisotropic surface energy per length, and 
$\omega$ is the orientation angle of the unit normal $\bm{n}$ exterior to the solid measured counterclockwise from the $x$ axis.
We choose an anisotropic surface energy function $\gamma(\omega)$ as the following classical model representative of the four-fold symmetry of a cubic crystalline structure (see, for example,
\cite{bertotti1998hysteresis}): 
\begin{equation} \label{eq:00}
\gamma(\omega)=\gamma_0\left(1+\varepsilon\cos{4\omega}\right), 
\end{equation}
where $\gamma_0$ is the mean interfacial tension, and the parameter $\varepsilon$ is a measure of the strength of anisotropy.  Thus
for $\varepsilon>0$, orientations with largest surface energy correspond to surface normals in the directions of the $\pm x$ and $\pm y$ axes with $\omega = 0, \pi/2, \pi, 3\pi/2$.  In the absence of stress the equilibrium shape is 4-fold symmetric, and, for sufficiently large anisotropy, the equilibrium shape has four corners in the $x$ and $y$ axis directions.  In the presence of biaxial stress with $\sigma^0_{xx} = \sigma^0_{yy}$, the 4-fold symmetry is preserved, however, if $\sigma^0_{xx} \ne \sigma^0_{yy}$ then the shape will only be 2-fold symmetric.   We allow for this asymmetry in the applied stress so that the effect of stress on the equilibrium shape can be more clearly resolved as a difference in behavior at corners in the $x$ and $y$ directions.  
Therefore, allowing for the 2-fold symmetric void shape we can describe the void shape in the first quadrant by
$r(\theta)$ on $0 \le \theta \le \pi/2$
and generate the shape in the other three quadrants by two-fold symmetry with respect to reflections across the $x$ and $y$ axes.

\subsection{Euler-Lagrange equations for minimum energy shape}
In \ref{Appendix:B} we derive the Euler-Lagrange equations for the shape that minimizes the total potential energy from the surface energy and elastic energy contributions,
\begin{equation}
\Pi=E_{\mbox{surf}}+\Pi_{\mbox{el}}  \label{eq:100}
\end{equation}
subject to the constraint of fixed void area $A$,
\begin{equation}
\int_0^{\pi/2} \frac{1}{2} r(\theta)^2 \, d\theta = A/4. \label{eq:area_constraint}
\end{equation}
The resulting equation for the shape is
\begin{equation}
\left(\gamma(\omega)+\gamma''(\omega) \right) \kappa - \frac{(1-\nu^2)}{2E} ( \sigma_{xx} + \sigma_{yy})^2 = \mu \hs \mbox{on $R=r(\theta)$}
\label{surface_condition}
\end{equation}
where $\kappa$ is the surface curvature given in \ref{Appendix:B}, and $\mu$ is the constant Lagrange multiplier corresponding to the surface chemical potential. 
The natural boundary conditions from the energy minimization give the corner angle conditions for $\omega$ at the ends of the interval 
(See Eq.\,\eqref{corner_bc} in \ref{Appendix:B}),
\begin{equation}
\frac{\gamma'(\omega)}{\gamma(\omega)} = \frac{\dot{r}(\theta)}{r(\theta)}  \hs \mbox{at $\theta=0, \pi/2$}, \label{eq:14}
\end{equation}
where $\dot{r}$ is $dr/d\theta$.

Note that the corner angle conditions are independent of stress, which is agreement with the conclusions of the scaling analysis of the energy in the neighborhood of the corner in Srolovitz and Davis \cite{srolovitz2001}.   
Note also that in the absence of stress Eq.\,\eqref{surface_condition} reduces to Herring's equation \cite{herring1951} with the appropriate corner angle conditions \cite{pimpinelli1998physics}. 

\subsection{Nondimensionalization}
We denote nondimensional variables with a tilde and use length scale $a = \sqrt{A/\pi}$ corresponding to the void size, stress scale $\sigma_{xx}^0$ from the applied stress, and strain scale $\sigma_{xx}^0 (1+\nu)/E$. 
The nondimensional variables are then defined by
\begin{eqnarray}
&r = a \tilde{r}, \hs \bm{\sigma} = \sigma_{xx}^0 \bm{\tilde{\sigma}}, \hs
\bm{\epsilon} = (\sigma_{xx}^0 (1+\nu)/E) \bm{\tilde{\epsilon}}, 
\hs \bm{u} = (a \sigma_{xx}^0 (1+\nu)/E) \bm{\tilde{u}} & \nonumber \\ 
&\kappa = (1/a) \tilde{\kappa}, \hs \gamma = \gamma_0 \tilde{\gamma}, \hs
\mu = (\gamma_0/a) \tilde{\mu}, \hs \Pi + C = (a \gamma_0) \tilde{\Pi}. &
\end{eqnarray}
Defining the nondimensional parameters
\begin{equation}
\chi = \frac{\sigma_{yy}^0}{\sigma_{xx}^0}, \hs 
\Lambda = \frac{2(\sigma_{xx}^0)^2 (1-\nu^2)a}{E \gamma_0},
\end{equation}
we rewrite Eqns.\,\eqref{eq:elasticity}-\eqref{eq:constitutive_law} in nondimensional form and drop the tilde notation to obtain the elasticity equations
\begin{gather}
\sigma_{ij,j} = 0 \hs \mbox{in $D$}, \\
\sigma_{ij} n_j = 0 \hs \mbox{on $\partial D$},\\
\sigma_{ij} \rightarrow \sigma^{\infty}_{ij} \hs \mbox{as $\sqrt{x^2 + y^2} \rightarrow \infty$}, \\
\sigma^{\infty}_{xx}=1, \hs \sigma^{\infty}_{yy} = \chi, \\
\epsilon_{ij} = \tfrac{1}{2} ( u_{i,j} + u_{j,i} ), \\
\epsilon_{11,22} - 2 \epsilon_{12,12} + \epsilon_{22,11}=0,\\
\sigma_{ij} =  \epsilon_{ij} + \frac{\nu}{1-2\nu} \delta_{ij}\epsilon_{kk}, 
\end{gather}
and the Eqns.\,\eqref{eq:area_constraint}-\eqref{eq:14} for the shape of the void become
\begin{gather}
\left(\gamma(\omega)+\gamma''(\omega) \right) \kappa - \tfrac{1}{4}\Lambda ( \sigma_{xx} + \sigma_{yy})^2 = \mu \hs \mbox{on $\partial D$}, \label{eq:nondim_surface}\\
\frac{\gamma'(\omega)}{\gamma(\omega)} = \frac{\dot{r}(\theta)}{r(\theta)}  \hs \mbox{at $\theta=0, \pi/2$}, \label{eq:nondimension_BC}\\
\int_0^{\pi/2} \frac{1}{2} r(\theta)^2 \, d\theta = \pi/4, \label{eq:nondimension_area}
\end{gather}
where the nondimensional total energy is
\begin{equation}
\Pi = \int_{\partial D} \gamma(\omega) \, ds + \frac{\Lambda}{4(1-\nu)} \int_{\partial D} n_i \sigma^{\infty}_{ij} u_j \, ds.
\label{eq:nondim_total_energy}
\end{equation}
All variables are nondimensional throughout the remainder of this paper.

\subsection{Boundary integral elasticity}

The elasticity problem can be formulated as a boundary integral equation for a single complex function (see e.g. \cite{mikhlin1957integral}) which is efficient for numerical implementation. The detailed derivation of the integral equation appears in Wang and Spencer \cite{wang2021}. Here we give an outline of this derivation.

First, introduce the Airy stress function $W(x,y)$ as a smooth function defined on $D$
and $\partial D$, such that $\sigma_{xx}=\partial^{2}W/\partial y^{2}$,
$\sigma_{xy}=-\partial^{2}W/\partial x\partial y$, $\sigma_{yy}=\partial^{2}W/\partial x^{2}$
(see \cite{gonzalez2008first}). Then $W(x,y)$ satisfies mechanical equilibrium in $D$, and the compatibility condition becomes the biharmonic equation:
\begin{equation}
\frac{\partial^{4}W}{\partial x^{4}}+2\frac{\partial^{4}W}{\partial x^{2}\partial y^{2}}+\frac{\partial^{4}W}{\partial y^{4}}=0.
\label{eq:biharmonic}
\end{equation}

Since $W(x,y)$ is a biharmonic function, following \cite{mikhlin1957integral}, we introduce the complex variable $z = x + i y$ and represent the biharmonic function $W(x,y)$ using two functions $\phi(z)$ and $\psi(z)$ (called Goursat
functions) which are holomorphic on $D$ and $\partial D$. Let $\text{\ensuremath{W}}(\text{\ensuremath{x}},\text{\ensuremath{y}})=\text{\mbox{Re}}\left\{ \bar{z}\phi(z)+\varsigma(z)\right\} $
and $\psi(z)=\varsigma'(z)$, where prime denotes a derivative with respect to $z$.
The nondimensional stress and displacement components can then be written as
\begin{gather}
\sigma_{xx}+\sigma_{yy}=4\mbox{Re}\left\{ \phi'(z)\right\} ,\label{eq:21}\\
\sigma_{yy}-\sigma_{xx}+2i\sigma_{xy}=2\left[\bar{z}\phi''(z)+\psi'(z)\right],\\
u_x + i u_y = 4(1-\nu) \phi(z).\label{eq:22}
\end{gather}
Since no external force is applied on $\partial D$, the boundary condition
on $\partial D$ is given by \cite{Muskhelishvili} as
\begin{equation}
\phi(z)+z\overline{\phi'(z)}+\overline{\psi(z)}=0\quad\mbox{\mbox{on}}\;z\in\partial D.\label{eq:26}
\end{equation}
We note that the far-field conditions for stress imply that 
$\phi'(z) \rightarrow (1 + \chi)/4$ and
$\psi'(z) \rightarrow (\chi -1)/2$ as 
$|z| \rightarrow \infty$.
We therefore let  
\begin{equation}
\phi(z)=(1+\chi)z/4+\varphi(z),\hs
\psi(z)=(\chi-1)z/2+h(z),
\end{equation} \label{eq:37}
where 
\begin{equation}
\varphi(z) \rightarrow 0,  \hs 
h(z) \rightarrow 0 \hs \mbox{as 
$|z| \rightarrow \infty$}, \label{eq:38}
\end{equation}
so that $\varphi$ and $h$ are analytic
on the region $D\cup\infty$.
Then Eq.\,\eqref{eq:26} on boundary is
\begin{equation}
    \overline{\varphi(z)}+\frac{1+\chi}{2}\overline{z}+\overline{z}\varphi'(z)+h(z)
+\frac{\chi-1}{2}z=0\quad\mbox{on}\;z\in\partial D.\label{eq:29}
\end{equation}
Multiply both sides of Eq.\,\eqref{eq:29}
by the factor $1/2\pi i\cdot \mathrm{d}z/(z-t)$, where $t$ is an arbitrary
point in the void region $D'$, and integrate along boundary $\partial D$, denoting
the integration contour $L$ as $\partial D$ traversed in the counterclockwise
direction. 
Eq.\,\eqref{eq:29} then becomes an integral equation which does not involve $h(z)$:
\begin{equation}
\frac{1}{2\pi i}\underset{L}{\int}\frac{\overline{\varphi(z)}}{z-t}\,\mathrm{d}z+\frac{1+\chi}{2\pi i\cdot2}\underset{L}{\int}\frac{\overline{z}}{z-t}\,\mathrm{d}z
+\frac{1}{2\pi i}\underset{L}{\int}\frac{\overline{z}\varphi'(z)}{z-t}\,\mathrm{d}z+\frac{\chi-1}{2}t=0.
\end{equation}
Now letting $t\rightarrow z_{0}$, where $z_{0}$ is a point on boundary
$\partial D$. Thus we obtain a singular
boundary integro-differential equation on $\partial D$:
\begin{equation}
\frac{1}{2}\overline{\varphi(z_{0})}+\frac{1}{2\pi i}\underset{L}{\int}\frac{\overline{\varphi(z)}}{z-z_{0}}\,\mathrm{d}z+\frac{1+\chi}{4}\overline{z_{0}}+\frac{1+\chi}{4\pi i}\underset{L}{\int}\frac{\overline{z}}{z-z_{0}}\,\mathrm{d}z
+\frac{1}{2}\overline{z_{0}}\varphi'(z_{0})+\frac{1}{2\pi i}\underset{L}{\int}\frac{\overline{z}\varphi'(z)}{z-z_{0}}\,\mathrm{d}z+\frac{\chi-1}{2}z_{0}=0.\label{eq:34}
\end{equation}
Once the integro-differential equation is solved for $\varphi(z)$ on
the boundary $\partial D$, Eq.\,\eqref{eq:29} determines $h$ on the boundary.
The stress at any point inside the solid can then be determined
by analytic continuation of boundary values $\varphi(z)$ and $h(z)$ into the domain $D$.

The free boundary problem Eq.\,\eqref{eq:nondim_surface} is now reduced to solving for the shape $r(\theta)$ that satisfies the surface equation
\begin{equation}
[\gamma(\omega)+\gamma''(\omega)]\kappa  -\frac{1}{4}\Lambda
[ 1 + \chi + 4 \mbox{Re}\{\varphi'(z)\} ]^2=\mu \quad \mbox{on}\ \partial D, \label{eq:114}
\end{equation}
subject to the boundary conditions \eqref{eq:nondimension_BC} and constraint \eqref{eq:nondimension_area}, where the 
total energy of the void from Eq.\,\eqref{eq:nondim_total_energy} is
\begin{equation}
\Pi = \int_{\partial D} \gamma(\omega) \, \mathrm{d}s 
+ \Lambda \int_{\partial D}
 n_1 \left[ \tfrac{1}{4}(1+\chi)x + \varphi_1 \right] +
 \chi n_2 \left[ \tfrac{1}{4}(1+\chi)y + \varphi_2 \right] \, \mathrm{d}s \label{total_energy}
\end{equation}
and where $\varphi_1 = \mbox{Re}(\varphi)$ and $\varphi_2 = \mbox{Im}(\varphi)$.

%% file: Sec3_Numerical_Method.tex
\section{Numerical Method}\label{Sec:3}
 
We develop a numerical method for simultaneously solving the coupled elasticity equations and free boundary problem based on the spectral collocation method with modifications to adapt to the singularity in the elastic stress at the corners of the void. 
The numerical method includes the following components: (i) an elasticity solver (see \cite{wang2021}) to solve $\varphi$ on the boundary
for a given shape; 
(ii) a calculation of the equilibrium shape satisfying Euler-Lagrange equation with given corner angles (which accesses (i) to compute $\varphi$
 to determine the elastic energy for candidate shapes); and (iii) numerical evaluation of the total energy for equilibrium shapes corresponding to 
different corner angles to find the corner angles minimizing the total energy (which accesses (ii) to determine the equilibrium shape for each choice of corner angle). 
Below we describe each component of the numerical method and summarize the validation of the numerical method with appropriate convergence tests.  

\subsection{Elasticity solver} \label{Sec:3.1}

This subsection briefly explains the elasticity solver solving the exterior plane-strain  elasticity problem for an arbitrary void shape $r(\theta)$. Complete details of this solver can be found in Wang and Spencer \cite{wang2021}. 

We use a Chebyshev basis to represent the real and imaginary parts
of the boundary values of the Goursat function on the first quadrant:
\begin{equation}
\varphi(\theta)=\sum_{k=1}^{N}a_{k}T_{k-1}(\theta)+i\cdot\sum_{k=1}^{N}b_{k}T_{k-1}(\theta) \qquad \theta\in [0,\pi/2],\label{eq:39}
\end{equation}
where $a_{k}$, $b_{k}$ are unknowns, and $T_{k}(\theta)$ are Chebyshev polynomials on $[0,\pi/2]$ obtained by applying the linear transform
\begin{equation}
    \theta=\frac{\pi}{4}+\frac{\pi}{2}x
\end{equation}
to standard Chebyshev polynomials $T_k(x)$ on $[-1,1]$. If the shape of the void has
corners, from the local asymptotic analysis near corners in \ref{Appendix:C} (or Sec.\,2.2 in \cite{wang2021}), the Goursat function $\varphi$ near the corner can not be well
approximated only by polynomials. Thus, considering the example where
the corner is located at $\theta_{c}=\pi/2$, we add a corner term
in $s^{\lambda-1}$, where $s$ is the arc length from the corner
and $\lambda$ is the strength of the elastic singularity given by the solution of Eq.\,\eqref{order_singularity}. 
We thus have
$\varphi(\theta)\sim c_{1}s^{\lambda-1}\sim c_{2}|\theta_{c}-\theta|^{\lambda-1}$
near the corner where $c_{1},c_{2}$ are constants.  To accommodate the 
local behavior near the corner we thus modify the expansion in Eq.\,\eqref{eq:39} as
\begin{equation}
\varphi(\theta)=a_{1}\theta^{\lambda_1-1}+a_{2}\left(\frac{\pi}{2}-\theta\right)^{\lambda_2-1}+\sum_{k=3}^{N}a_{k}T_{k-3}(\theta)\\
+i\cdot\left(b_{1}\theta^{\lambda_1-1}+b_{2}\left(\frac{\pi}{2}-\theta\right)^{\lambda_2-1}+\sum_{k=3}^{N}b_{k}T_{k-3}(\theta)\right), \label{eq:115}
\end{equation}
where $\lambda_1$ and $\lambda_2$ are the singularity exponents at $\theta=0$ and $\theta=\pi/2$, respectively.

Since $\varphi$ is a continuous function of $\theta$, using the symmetry of $\varphi$, we have the
continuity conditions on both ends of the interval
\begin{gather}
\mbox{Re}\left\{ \varphi(\theta)\right\} =0\quad\mbox{at}\;\theta=\pi/2,\label{eq:42}\\
\mbox{Im}\left\{ \varphi(\theta)\right\} =0\quad\mbox{at}\;\theta=0\text{.}\label{eq:43}
\end{gather}

The boundary integral equation \eqref{eq:34} holds for any $z_{0}$ on $L$. Following the idea of the spectral collocation method, we pick $z_{0}=z_{i}=z(\theta_{i})$, 
$i=1,2,\ldots,N-1$ as collocation points, where $\theta_{i}=\pi(x_{i}+1)/4$
and $x_{i}$ is the i-th root of degree ($N-1$) Legendre polynomial
$P_{N-1}$ (Gauss-Legendre quadrature points). The choice of the number
of collocation points is determined by the number of unknowns. We enforce
both real part and imaginary part of Eq.\,\eqref{eq:34} at each $z_{i}$
($2N-2$ equations) together with two boundary conditions \eqref{eq:42} and
\eqref{eq:43} to give a system of $2N$ equations for $2N$
unknown coefficients $a_{k}$, $b_{k}$. 

The singularity of the Cauchy principal value integrals in Eq.\,\eqref{eq:34}
is extracted using
\begin{equation}
\underset{L}{\int}\frac{\overline{\varphi(z)}}{z-z_{i}}\,\mathrm{d}z  =\underset{L}{\int}\frac{\overline{\varphi(z_{i})}}{z-z_{i}}\,\mathrm{d}z+\underset{L}{\int}\frac{\overline{\varphi(z)}-\overline{\varphi(z_{i})}}{z-z_{i}}\,\mathrm{d}z
 = \pi i\thinspace\overline{\varphi(z_{i})}
+\underset{L}{\int}\frac{\overline{\varphi(z)}-\overline{\varphi(z_{i})}}{z-z_{i}}\,\mathrm{d}z. \label{eq:44}
\end{equation}
All integrals in Eq.\,\eqref{eq:34}
are then evaluated numerically using nested Gauss-Legendre quadrature (see \cite{wang2021,Hoskins2019} for details). $\varphi(\theta_{i})$ can
be obtained from the Chebyshev coefficients by evaluating the Chebyshev series
at the collocation points $\theta=\theta_{i}$. 

In solving the system of equations for the Chebyshev coefficients, we apply the following analyticity constraint: 
\begin{equation}
\frac{1}{2\pi i}\underset{L}{\int}\frac{\varphi(z)}{z-z_{i}}\thinspace \mathrm{d}z+\frac{1}{2}\varphi(z_{i})=0\label{eq:50}
\end{equation}
for any collocation point $z_{i}$. This constraint is necessary to eliminate numerical solutions $\varphi(\theta)$ which satisfy the integral equation \eqref{eq:34} but are not analytic in $D$ (see \cite{wang2021} for details). Then Eq.\,\eqref{eq:50}
at collocation points is a set of linear equations of $a_{k}$ and
$b_{k}$. We therefore
use Eq.\,\eqref{eq:50} at each collocation point in addition to
our boundary integral equation \eqref{eq:34} to construct analytic $\varphi(z)$. Finally, by combining the $2(N-1)$ equations from Eq.$\,$\eqref{eq:50} at
collocation points, the $2(N-1)$ equations from discretized integral
equation Eq.$\,$\eqref{eq:34}, and two equations Eqs.$\,$\eqref{eq:42},\,\eqref{eq:43} from
boundary conditions at each end of the domain, we obtain an overdetermined
linear system of $(4N-2)$ equations for the $2N$ unknowns $a_{k}$
and $b_{k}$. 
More discussions of analyticity constraints Eq.\,\eqref{eq:50} and the uniqueness of the solution from the elasticity solver can be found in \cite{wang2021}.

\subsection{Solving Euler-Lagrange equations for fixed corner angles} \label{Sec:3.2}

We use the spectral collocation method \cite{hesthaven2007spectral} to determine the equilibrium shape assuming the solid corner angles $\alpha_1, \alpha_2$ (see Fig.\,\ref{fig:1}) are given.
We represent the surface $r(\theta)$ in polar coordinates by a set of Chebyshev polynomials with corner terms as in the expansion of $\varphi$ given in Eq.\,\eqref{eq:115}. In \ref{Appendix:D} we use an asymptotic analysis of the surface equilibrium equation \eqref{eq:114} to determine the functional form of the shape in response to the singular elastic stress at the corner as
\begin{equation}
r(\theta)=c_1\theta^{2\lambda_1-2}+c_2\theta^{\lambda_1}+c_3\left(\frac{\pi}{2}-\theta\right)^{2\lambda_2-2}+
c_4\left(\frac{\pi}{2}-\theta\right)^{\lambda_2}+\sum_{k=5}^{N}c_{k}T_{k-5}(\theta), \\ \label{eq:113}
\end{equation}
where $c_1,c_2,c_3,c_4$ are determined from the elastic singularity as given in \ref{Appendix:D}, $c_5,\ldots,c_N$ are unknowns, $T_{k}(\theta)$ are Chebyshev polynomials, and collocation points $\theta_i$ are the $(N-5)$ Legendre points similar to Sec.\,\ref{Sec:3.1}. Note the Legendre points are more dense near the ends of the interval which is beneficial to analyzing the corner behavior since the corners occur at $\theta=0,\ \pi/2$. 

Applying the equilibrium condition Eq.\,\eqref{eq:114} at $\theta_i$ we obtain the nodal equations for the shape
\begin{equation}
    \left[1-15\varepsilon\cos{4\left(\arctan{\frac{\dot{r_i}\sin{\theta_i}+r\cos{\theta_i}}{\dot{r_i}\cos{\theta_i}-r\sin{\theta_i}}}\right)}\right]
    \frac{2\dot{r_i}^2-r_i\ddot{r_i}+r_i^2}{(\dot{r_i}^2+r_i^2)^{3/2}}-\frac{1}{4}\Lambda(4\,\mbox{Re}\left\{\frac{\dot{\varphi_i}}{\dot{z_i}}\right\}+(1+\chi) )^2=\mu,
\end{equation}
where $r_i,\ \dot{r}_i,\ \ddot{r}_i,\ \dot{\varphi}_i,\ \dot{z}_i$ are $r,\ \dot{r},\ \ddot{r},\ \dot{\varphi},\ \dot{z}$ at $\theta=\theta_i$. Here $z_i=r_i\cos(\theta_i)+\sqrt{-1}\ r_i \sin(\theta_i)$, and $\varphi$ is determined from the elasticity solver for a given the shape $r(\theta)$. 

The two boundary conditions for $r(\theta)$ are the imposed solid corner angles $\alpha_1$ and $\alpha_2$ at $\theta=0,\,\pi/2$ respectively determined by geometry as
\begin{equation}
    \frac{r'(0)}{r(0)}=-\tan\left(\frac{\alpha_1-\pi}{2}\right),\qquad
    \frac{r'(\pi/2)}{r(\pi/2)}=\tan\left(\frac{\alpha_2-\pi}{2}\right).
\end{equation}
After rewriting these equations (4 equations from corner terms, $(N-5)$ equations at collocation points, 2 equations from boundary conditions) in terms of the $c_i$ and $\mu$, it can be seen that the system of equations is a nonlinear system
with $(N+1)$ unknowns and $(N+1)$ equations. We can determine the unknown coefficients $c_i$ and $\mu$ by solving the nonlinear system using MATLAB \emph{fsolve} \cite{dennis1981adaptive,MatlabOTB}.

\subsection{Total energy}

We now consider the total potential energy as a function of the imposed corner angles to confirm the corner angle conditions from the Euler-Lagrange equations Eq.\,\eqref{eq:14} and to validate our numerical method. \ref{Appendix:A} states that the minimization of total potential energy is equivalent to minimizing the expression in Eq.\,\eqref{total_energy}. By using the numerical results of the equilibrium shape and $\varphi$ for fixed corner angles, we can evaluate the integrals in Eq.\,\eqref{total_energy} using numerical quadrature. Then, we can compare the total energy of equilibrium shape for different corner angles. By the principle of minimum energy, the correct corner angles give minimum energy.
The numerical solution of the equilibrium corner angles is therefore
\begin{equation}
    \alpha^*_1,\alpha^*_2=\underset{\alpha_1,\alpha_2}{\mbox{arg\,min}} \,\left\{\Pi(\alpha_1,\alpha_2)\right\} \label{eq:18}
\end{equation}
where $\Pi(\alpha_1,\alpha_2)$ is the total potential energy of the equilibrium shape corresponding to corner angles $\alpha_1,\,\alpha_2$. We show in detail later that these numerically determined corner angles are 
independent of the elastic singularity and are the same as the corner angles prescribed by the Euler-Lagrange equation Eq.\,\eqref{eq:14}.

\subsection{Validation of numerical methods}
Before using the numerical method to find the elasticity solution, free boundary shape and corner angles, we validate the numerical method in the following examples. 
\subsubsection{Convergence to Wulff shape (without stress) for known corner angles}

Wulff's shape \cite{wulff1901} is the exact solution for the void shape in the absence of stress. Since no stress is applied, the corner angles at $\theta=0$ and $\theta=\pi/2$ are equal ($\alpha=\alpha_1=\alpha_2$). The corner angle for the Wulff shape is the solution $\alpha$ to the following equation  \cite{cabrera1964}:
\begin{equation}
\tan\left((\alpha-\pi)/2\right)=-\frac{\gamma'((\alpha-\pi)/2)}{\gamma((\alpha-\pi)/2)}. \label{eq:17}
\end{equation}
With the exact corner angle $\alpha$ from Eq.\,\eqref{eq:17}, we set corner angles in the boundary conditions Eq.\,\eqref{eq:14} as known ($\omega=\alpha_1/2-\pi/2=\alpha_2/2-\pi/2$) and apply the
numerical method to find the equilibrium shape of the surface. This experiment tests the accuracy of finding the shape with fixed corner angles.
\begin{figure}[!t]
    \centering
    \captionsetup[subfigure]{justification=centering}
    \begin{subfigure}{.5\textwidth}
        \centering
        \includegraphics[width=75mm]{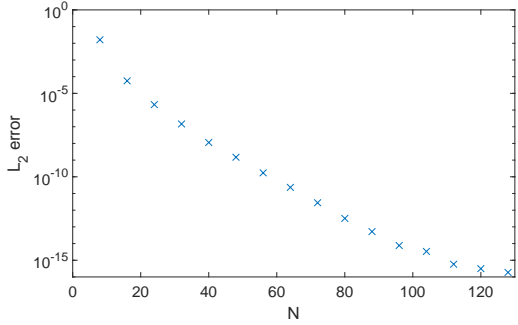}
        \caption{$L_2$ error vs. number of collocation points $N$.} \label{fig:8}
    \end{subfigure}%
    \begin{subfigure}{.5\textwidth}
        \centering
        \includegraphics[width=75mm]{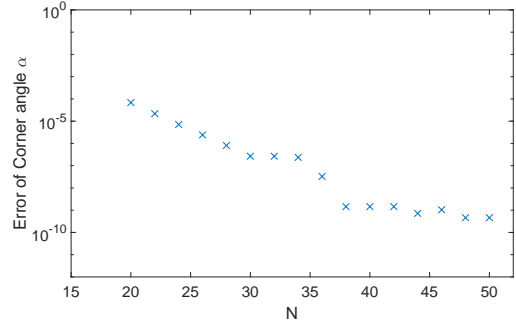}
        \caption{Error of corner angle vs. number of collocation points $N$.}
        \label{fig:9}
    \end{subfigure}
    \caption{Convergence to Wulff shape (without stress), $\varepsilon=0.08$.}
\end{figure}

Fig.\,\ref{fig:8} compares the numerical solution with the exact solution given by the Wulff shape and shows that the numerical solution converges to the exact solution with the rate of convergence close to exponential convergence. Here the $L_2$ error is defined by
\begin{equation}
L_2\ \mbox{error}=\left[\frac{2}{\pi}\int_{0}^{\pi/2}\Big|r(\theta)-r_{\mbox{exact}}(\theta)\Big|^{2}\,\mathrm{d}\theta\right]^{1/2}. \label{eq:20}
\end{equation}
We will use this definition of $L_2$ error in the rest of our paper.

\subsubsection{Convergence of corner angle from energy minimization (without stress)}

Now assume that the corner angles of the equilibrium shape are unknown. We find the corner angles by minimizing the total potential energy in Eq.\,\eqref{eq:18}
and compare the solution with the exact corner angle in Eq.\,\eqref{eq:17}. 
The corner angle determined by our numerical method converges to the exact corner angle of the Wulff shape as shown in Fig.\,\ref{fig:9}.
Thus, Fig.\,\ref{fig:9} confirms that the exact corner angle gives the minimum
surface energy of all corner angles for the void shape without stress.

\subsubsection{Convergence of elasticity solver (fixed shape with corners)}

\begin{figure}[!t]
   \centering
   \includegraphics[width=75mm]{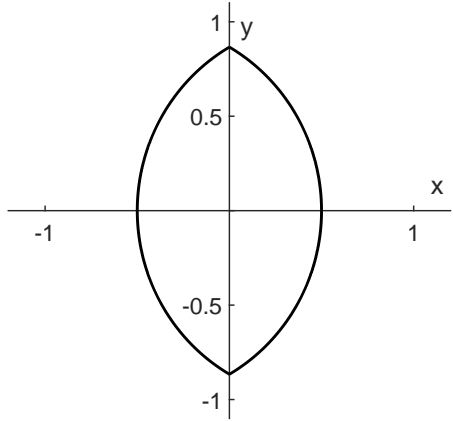}
   \caption{\label{fig:23}Overlapping circles shape   ($\alpha_0=2\pi/3$).}
\end{figure}

We now present the convergence of the elasticity solver using a void shape with
a void corner angle less than $\pi$ as shown in Fig.\,\ref{fig:23}. The void shape is given by taking the $x>0$ portion of the unit circle with center
at $(\cos(\alpha_0),0)$ and reflecting it across the $y$ axis. Thus
the parameter $\alpha_0$ controls the amount of overlap between the
circles. The equation of the overlapping circles void shape in polar
coordinates is $r(\theta)=\cos\alpha_0\cos\theta+\sqrt{1-\sin^{2}\theta\cos^{2}\alpha_0}$
in the first quadrant. We consider the case of uniaxial applied stress 
$\chi=0$.  Due to the corner at $\theta=\pm\pi/2$ there is an elastic singularity at the corner. The exact solution for $\sigma_{xx}+\sigma_{yy}$
is given in \cite{ling1948stresses}. Fig.\,\ref{fig:29} shows the $L^{2}$ error of $\sigma_{xx}+\sigma_{yy}$ on the surface versus number
of collocation points. The absolute $L^{2}$  error is less than $10^{-4}$ when $N$ is large.
So, even in this case with a stress singularity
at the corner, the numerical solution is of good accuracy for moderately large $N$. See \cite{wang2021} for a more detailed testing of the elasticity solver.

\subsubsection{Convergence of shape with stress for fixed corner angle}

\begin{figure}[!t]
   \centering
   \begin{minipage}{.5\textwidth}
       \includegraphics[width=75mm]{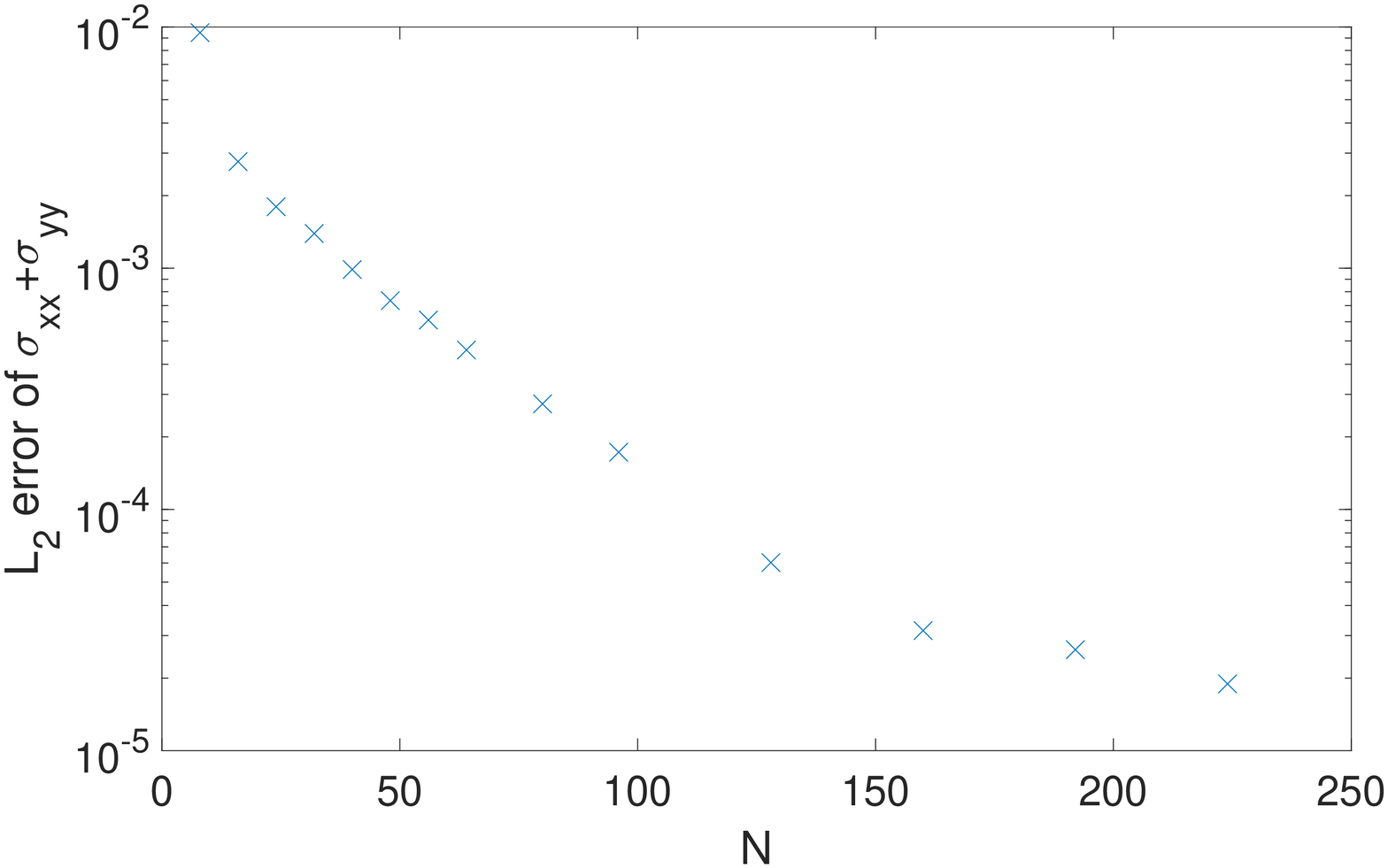}
       \caption{\label{fig:29}Overlapping circle case ($\alpha_0=2\pi/3$): $L^{2}$ error of \\\hspace{\textwidth} $\sigma_{xx}+\sigma_{yy}$ versus number of collocation points $N$ for $\chi=0$.}
    \end{minipage}%
    \begin{minipage}{.5\textwidth}
        \centering
        \includegraphics[width=65mm]{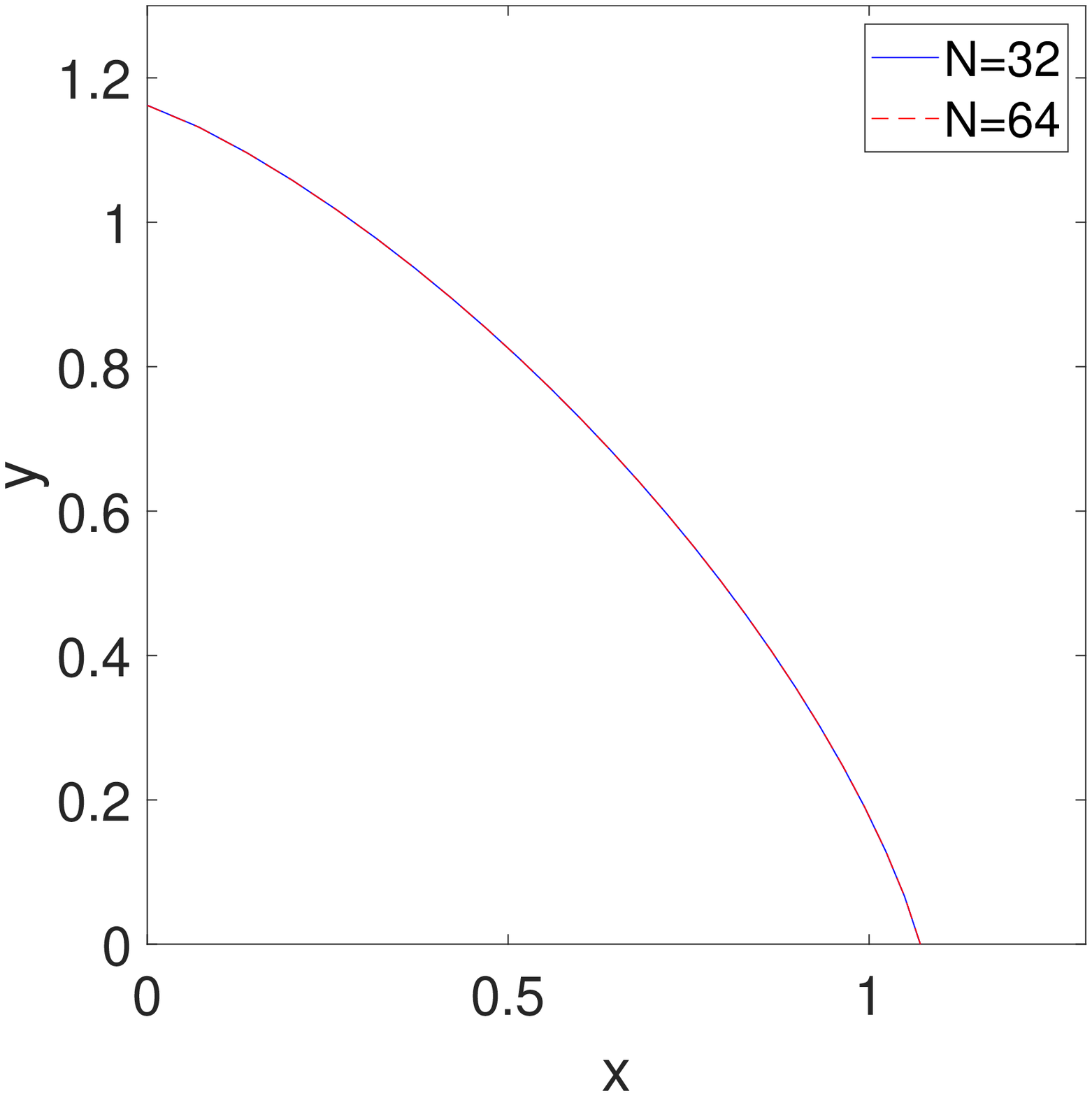}
     \caption[Convergence of shape.]{\label{fig:33} Calculated void shape (in first quadrant), $N=32$ (solid blue) and $N=64$ (red dashed), for $\varepsilon=0.08,\ \chi=0,\ \Lambda=0.15$.  The total energies are 
     5.852825050913609 ($N=32$) and  
     5.852823956354603 ($N=64$) which agree to 6 digits.}
    \end{minipage}
\end{figure}

We test the convergence of the shape when stress is applied at infinity and the corner angles are fixed. The convergence of the shape is given in Fig.\,\ref{fig:33} where the shapes for $N=32$ and $N=64$ are indistinguishable.  Also given in the caption are the numerically calculated total energies for $N=32$ and $N=64$ which agree to 6 digits.  

%% file: Sec4_Result.tex
\section{Corner Angle Problem with Elasticity}\label{Sec:4}
 
Now that we have validated our numerical method to confirm that (i) in the absence of stress we recover the equilibrium void shape and the corner angles,
(ii) our elasticity solver accurately resolves stresses for the singular elastic fields associated with a corner, and (iii) the void shape and energy in the presence of elastic stress can be found for fixed corner angles, we investigate how the shape and energy for the stressed elastic void depends on corner angles. 

\subsection{Example results: convergence of minimum total energy vs. corner angles}

In this example, we fix the anisotropy of the surface energy $\varepsilon=0.08$, relative strength of the stress $\Lambda=0.15$ and $\chi=0$ as a prototype problem corresponding to uniaxial loading. 

\begin{figure}[!t]
    \centering
    \begin{minipage}{.5\textwidth}
    \includegraphics[width=75mm]{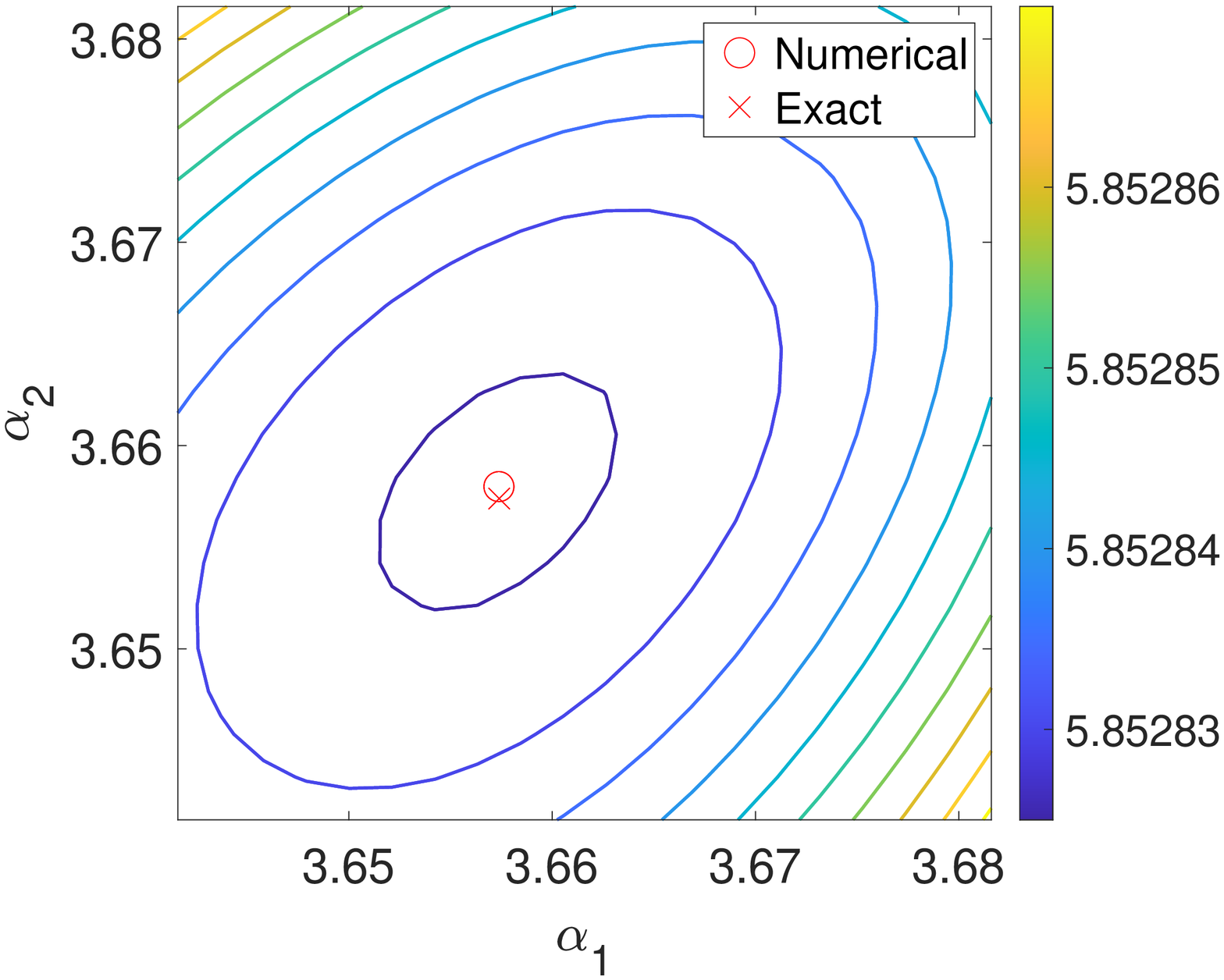}
    \caption[Level curve of minimizing total energy]{\label{fig:35}Level curve plot of the total energy.  $\varepsilon=0.08$,\\\hspace{\textwidth}
    $\chi=0,\ \Lambda=0.15,\ N=64$.}
    \end{minipage}%
    \begin{minipage}{.5\textwidth}
       \centering
       \includegraphics[width=75mm]{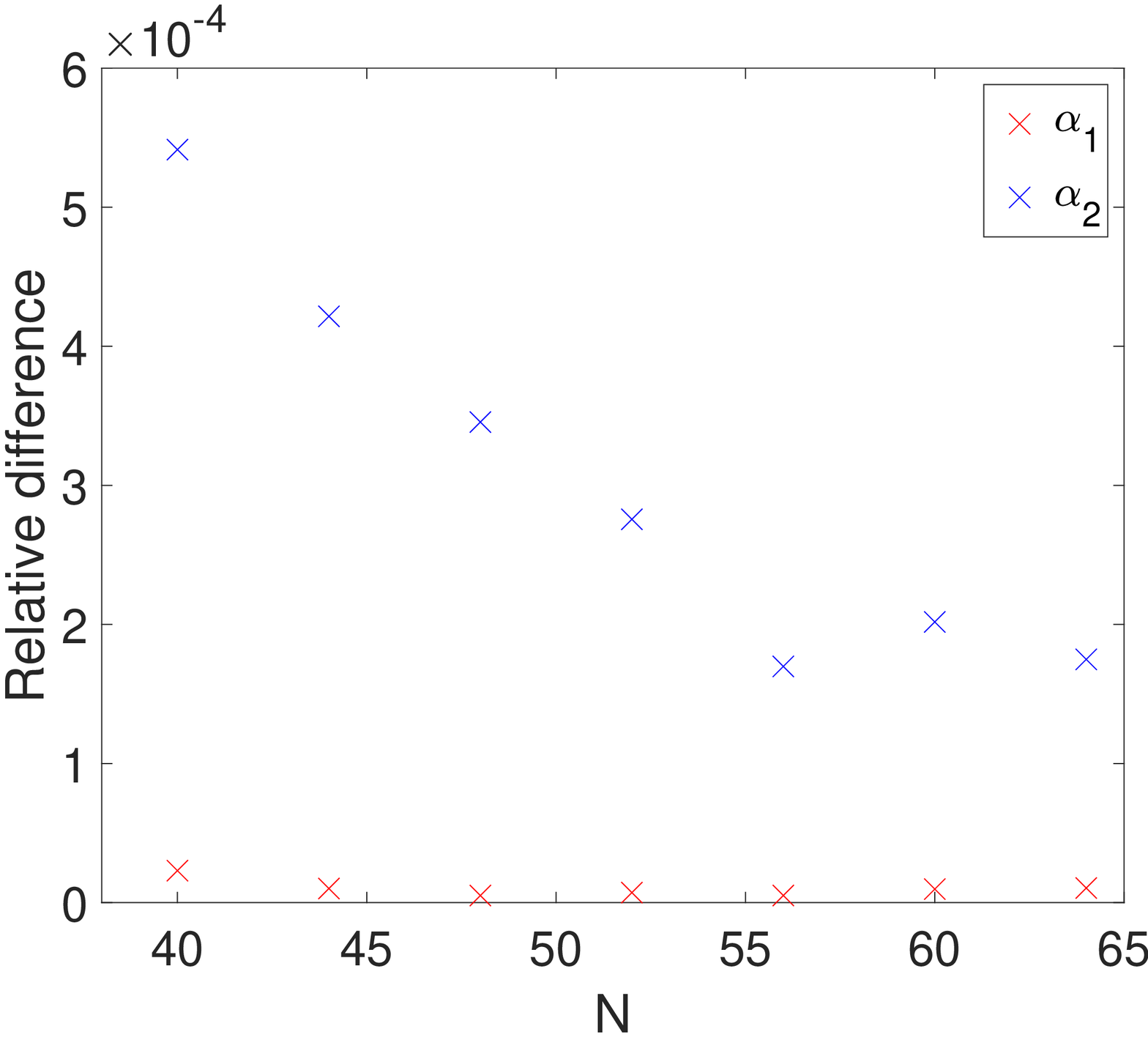}
       \caption[Relative difference vs. $N$]{\label{fig:34}Relative difference of corner angle from unstressed corner angle $\alpha_0$ vs. number of collocation points $N$,}
    \end{minipage}
\end{figure}

Fig.\,\ref{fig:35} is a level curve plot showing the total potential energy as a function of the imposed corner angles $\alpha_1$ and $\alpha_2$ evaluated on a grid in $\alpha_1,\alpha_2$ space. 
This numerical result shows the total energy is a smooth convex function that is concave up and has a minimum around $(3.66,3.66)$. 

Since we demonstrated the convergence of the numerical method with respect to the corner angles
and total energy, we can find the equilibrium corner angles by using our numerical method as part of the minimization of Eq.\,\eqref{total_energy}. Fig.\,\ref{fig:34} shows the relative difference between the calculated corner angles and the exact no-stress corner angles when $\Lambda=0.15$ as a function of $N$. The computed corner angles $\alpha_1,\,\alpha_2$  have a relative difference of less than 0.0005 compared to the non-stressed corner angle $\alpha_0$. 
This result illustrates that corner angles do not change from the no-stress value for a finite stress parameter $\Lambda$. The level curve plot in Fig.\,\ref{fig:35} shows the solution of the corner angle that minimizes the total energy.
The location of the minimum energy is nearly identical that for the non-stressed corner angle, differing by about $0.03\%$.  Since the magnitude of the elasticity terms is non-negligible (and in fact singular at the corner), 
we can thus conclude that stresses do not modify the equilibrium corner angles of the void.

\subsection{Analysis of the void shape} \label{Sec:5.6}

\begin{figure}[!t]
\centering
\includegraphics[width=100mm]{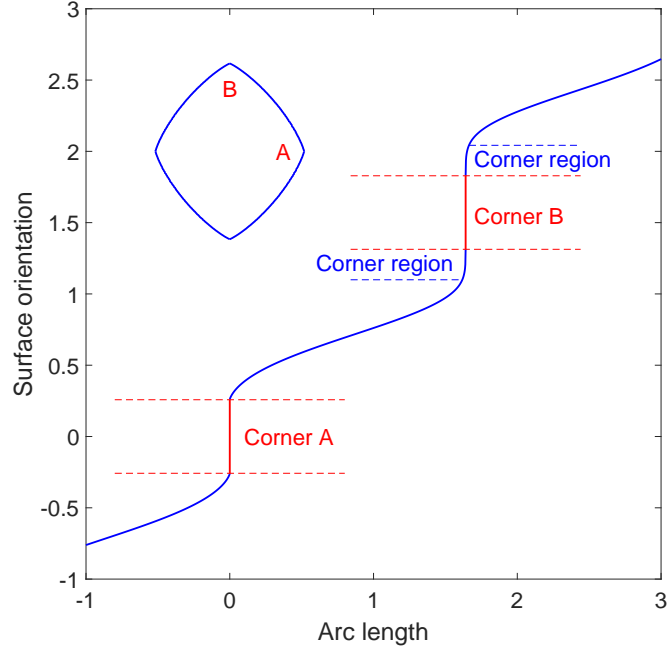}
\caption[Surface orientation vs. arc length]{\label{fig:37}Surface orientation $\omega$ vs. arc length $s$ (equilibrium shape in the top left of the figure), $\varepsilon=0.08,\ \chi=0,\ \Lambda=0.3$.}
\end{figure}

In this section we consider the numerically-determined equilibrium void shape and its behavior near the corners to explain the apparent discrepancy
between the regularized corner results in Siegel et al.\,\cite{siegel2004} and the sharp corner prediction of Srolovitz and Davis\,\cite{srolovitz2001}.
Fig.\,\ref{fig:37} (top left inset) shows the equilibrium shape for $\varepsilon=0.08,\ \chi=0,\ \Lambda=0.3$. From the figure, it is clear that the shape is elongated vertically and appears to have
a smaller corner angle at the top and the bottom of the void relative to the left and right sides of the void. 
This shape and asymmetric corner behavior is consistent with the result shown in Fig.\,7 
of Siegel et al.\,\cite{siegel2004}. But as has been shown in the previous section, the numerical values of the corner angles at all four void 
corners are identical, in agreement with the predictions of Srolovitz and Davis \cite{srolovitz2001} that the corner angles are independent of elastic stress.
 
To resolve this discrepancy, we examine the corner region more closely. Fig.\,\ref{fig:37} is a plot of the surface
orientation $w$ vs. arc length $s$ for the shape in Fig.\,\ref{fig:37}. We take the point with the maximum $x$
value as the starting point for the arc length, and the positive $x$-direction is the reference direction for the orientation angle. 
We analyze the direction change near each corner. The surface orientation
has a jump at each corner of the void (this follows from the definition of a corner). Corner A is the corner on the right side of the void. Corner B is the corner on the top of the void.
Both corners have the same jump in orientation (labeled in red) which is precisely the jump in orientation prescribed by the no-stress corner angles.
But a comparison of the behavior near the jump in orientation shows a distinct difference between the two corners.
In particular, the corner B has a rapid variation
of the orientation near the corner and includes a near-vertical portion connecting to the corner. This rapid variation in the orientation angle gives 
a solid corner angle that is effectively larger (void angle smaller) as shown on the equilibrium shape inset in Fig.\,\ref{fig:37}.

The difference in apparent corner angles and the variation of the orientation angle near the corner can be explained by the elastic energy. 
For corner A, since the applied stress is in the $x$-direction only, there is no elastic singularity at the corner and only a mild alteration
of the shape due to elasticity. However, at corner B the elastic stress is singular. By the asymptotic analysis
in \ref{Appendix:C}, the corner singularity of the elastic stress in the corner region leads to a matching singularity in the surface curvature $r''$. Thus the surface orientation, 
which depends on $r'$, is changing rapidly near the corner as in the blue-labeled corner regions in Fig.\,\ref{fig:37}.
Observing corner B on a large scale, corner B is a combination of the exact corner and two regions with singular curvature, which gives the appearance of as a corner with a larger corner angle than the exact value.

%% file: Sec5_Discussion.tex
\section{Discussion}
 
The results in Sec.~\ref{Sec:5.6} demonstrate that the elastic singularity associated with a corner results in an equilibrium shape with a singularity in the curvature near the corner that causes the apparent corner angle to be different from the corner angle associated with the stress-free Wulff shape.
The magnitude of this effect on the corner region depends on the nondimensional parameter 
\begin{equation}
\Lambda = \frac{2 (\sigma_{xx}^0)^2 (1-\nu^2) a }{E \gamma_0}.
\label{Lambda_Sec5}
\end{equation}
For the case of $\Lambda =0$, there is no stress effect and the shape is identical to the Wulff shape.  For $\Lambda = 0.3$ as in Fig.\ \ref{fig:37} the effect of elasticity on the shape and the deviation of the apparent corner angle from the stress-free Wulff-shape corner angle is clearly distinguishable.  

\begin{figure}[t]
\centering
\includegraphics[width=75mm]{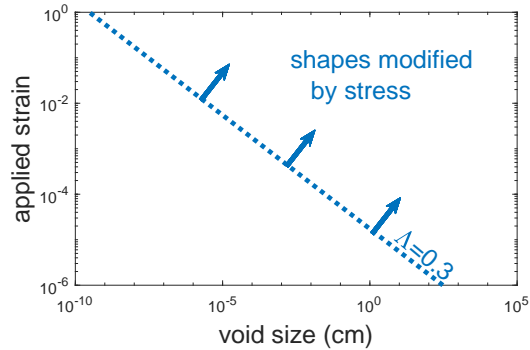}
\caption[Void size for which elastic effects affect the corner angle.]{\label{fig:void_size} Schematic of applied strains $\epsilon^0_{xx}$ and void sizes $a$ for which the void shape and corner angle are clearly modified by stress.  Here $l_0 = 2 \times 10^{-9} \mbox{cm}$. }
\end{figure}

From the dependence of $\Lambda$ on the dimensional parameters we can thus infer for fixed applied stress $\sigma^0_{xx}$ that
(i) in the limit of very small voids ($a \rightarrow 0$) we have $\Lambda \rightarrow 0$ and so small voids have shapes that are unaffected by 
elasticity (i.e. the usual case of surface effects dominating bulk effects at small length scales); and 
(ii) for large voids the effect of elasticity on the equilibrium void shape becomes important.  An estimate of the void size for which elasticity generates a noticeable change to the apparent corner angle (as in Fig.\ \ref{fig:37}) can be found by solving for the void size $a$ from Eq.\ \eqref{Lambda_Sec5} and rewriting the applied stress in terms of an applied strain $\epsilon^0_{xx}$ using the plane strain result
\begin{equation}
\sigma^0_{xx} = E \epsilon^0_{xx} / (1-\nu^2)
\end{equation}
to obtain
\begin{equation}
a = \frac{ \Lambda l_0}{2 (\epsilon^0_{xx})^2}
\end{equation}
where 
\begin{equation}
l_0 = \gamma_0 (1-\nu^2)/E
\end{equation}
is a material length scale associated with the ratio of surface energy to the elastic modulus.   The material length scale $l_0$ is roughly the length scale of the atomic spacing on the crystal lattice (see for example \cite{weir2008implications}). 
Taking a typical value for $l_0$ as $l_0 \approx 2 \times 10^{-9} \mbox{cm}$
and using $\Lambda = 0.3$ (corresponding to a noticeable change in the apparent corner angle as in Fig.\ \ref{fig:37}) we find that for an applied strain of $\epsilon^0_{xx} = 0.001$, voids of size $a=3 \mu\mbox{m}$ and larger would have shapes being significantly affected in the corner regions by elastic effects.   
Generalizing this idea, Fig.\ \ref{fig:void_size} is a schematic diagram that illustrates the regions in strain $\epsilon_{xx}^0$ and void size $a$ that exceed the nominal threshold of $\Lambda =0.3$ for which the effects of elasticity on the shape and corner angle are clearly visible.

In theory it would be useful to extract the magnitude of the apparent change to the corner angle from the form of the curvature singularity described in \ref{Appendix:D}. In particular, the singular contributions to the curvature $\kappa$ could be substituted into the Frenet formula $d \omega / ds = \kappa$ relating the curvature to the orientation angle and integrated to obtain the resultant change on $\omega$ in the corner region due to the singular terms in $\kappa$.  However, while the singular curvature terms do give a finite contribution to the orientation angle near the corner, there is no well-defined ``cut-off'' away from the corner for isolating the corner contributions because the singular terms in $\kappa$ 
decay slowly away from the corner region.   Thus, the best we can do in terms of quantifying the effect on the corner region is to use the singular terms in conjunction with all the other terms in the ``global'' expansion of the shape as part of our numerical solution, as shown in Fig.\ \ref{fig:37}.



%% file: Sec6_Conclusion.tex
\section{Conclusion} \label{Sec:6}

This work presents a numerical method to solve the free boundary problem for the energy-minimizing shape of a void inside a two-dimensional elastic solid with anisotropic surface energy.  In particular, we are interested in the effect of elasticity on corner angles of the void and address the apparent contradiction in the literature on whether elasticity affects the equilibrium void corner angle \cite{srolovitz2001,siegel2004}.

In the absence of elasticity we recover the classic Wulff shape.
In the presence of elasticity the corners generate singularities in the elastic field.  Using asymptotic solutions for the elastic singularity 
we designed a numerical method to solve the elasticity problem for a general class of void shapes with corners.
The elasticity solver was used as a component in the numerical solution of the free boundary problem for the void shape.  Asymptotic analysis of the corner region shows that the singularity in the elastic field generates singular terms in the curvature at the corner which we incorporate into our numerical solution.

The numerical results show that the precise corner angle of a void is unaffected by elasticity, in agreement with the scaling argument in Srolovitz and Davis \cite{srolovitz2001}, however, the rapid change in orientation near the corner due to the curvature singularity effectively changes the apparent corner angle when viewed macroscopically.  In this regard our results agree with the regularized corner results in Siegel et al \cite{siegel2004}.

Thus, our work resolves the apparent contradiction in Srolovitz and Davis \cite{srolovitz2001} and Siegel et al \cite{siegel2004}.  In addition to confirming that the true corner angle is independent of stress, a broad 
implication of our work is that determination of the equilibrium free boundary configuration of a stressed solid with inward pointing (notch-like) corners involves a nontrivial singular shape in the vicinity of the corner region and effectively modifies the macroscopic corner angle.

%% file: Appendix_A.tex
\section{Derivation of elastic potential energy} \label{Appendix:A}

Here we derive a simplified expression for the elastic potential energy of the elastically stressed void as given in Sec.~\ref{Sec:2}, Eqs.~\eqref{eq:elastic_potential_energy}-\eqref{eq:elastic_work} as
\begin{equation}
    \Pi_{\mbox{el}} = E_{\mbox{int}} - W_{\mbox{ext}} 
\end{equation}
where
\begin{equation}
E_{\mbox{int}} = \int \!\int_D \tfrac{1}{2}  \sigma_{ij} \epsilon_{ij} \, dA
\end{equation}
and
\begin{equation}
W_{\mbox{ext}} = \int_{L_\infty} n_i \sigma_{ij} u_j \, ds.
\end{equation}
Note that because of the unbounded domain, both $E_{\mbox{int}}$ and $W_{\mbox{ext}}$ are unbounded integrals and so it is necessary to simplify the result 
to eliminate taking the difference of unbounded energies. To this end we rewrite the integral in $E_{\mbox{int}}$ as
\begin{equation}
\begin{aligned}
E_{\mbox{int}} & = \frac{1}{2}\underset{D}{\iint}\epsilon_{ij}\sigma_{ij}\,\mathrm{d}A  =\frac{1}{4}\underset{D}{\iint}(u_{i,j}+u_{j,i})\sigma_{ij}\,\mathrm{d}A \\
& = \frac{1}{2}\underset{D}{\iint}u_{i,j}\sigma_{ij}\,\mathrm{d}A= \frac{1}{2}\underset{D}{\iint}\partial_j(u_i\sigma_{ij})-(\partial_j\sigma_{ij})u_i\,\mathrm{d}A \\
& =  \frac{1}{2}\underset{L_1+L_2}{\int}\sigma_{ij}u_i n_j\,\mathrm{d}s=  \frac{1}{2}\underset{L_2}{\int}\sigma_{ij}u_i n_j\,\mathrm{d}s  =  \frac{1}{2}\sigma_{ij}^\infty\underset{L_2}{\int}u_i n_j\,\mathrm{d}s \\
& = \frac{1}{2}\int_{L_\infty} n_i \sigma_{ij} u_j \, ds
=\frac{1}{2}W_{\mbox{ext}}\label{eq:103}
\end{aligned}
\end{equation}
where we applied the divergence theorem and the force balance condition Eq.\,\eqref{eq:3} on the third line, and where $L_1=\partial D$ is the interface between the void and the solid, and $L_2=L_{\infty}$ is the boundary at infinity. 
 
Following the approach in Appendix A of \cite{suo1994}, 
\begin{equation}
E_{\mbox{int}}=\frac{1}{2}\sigma_{ij}^\infty\underset{L_2}{\int}u_i n_j\,\mathrm{d}s=\frac{1}{2}\sigma_{ij}^\infty\underset{L_1+L_2}{\int}u_i n_j\,\mathrm{d}s-
\frac{1}{2}\sigma_{ij}^\infty\underset{L_1}{\int}u_i n_j\,\mathrm{d}s.
\label{eq:E_suo}
\end{equation}
The first integral on the right of Eq.\,\eqref{eq:E_suo} can be analyzed by assuming an auxiliary solid with the same geometry as the original solid. The auxiliary solid is under constant
stress $\sigma^{\infty}_{ij}$ everywhere. Then the first integral is the virtual work done by the traction on the auxiliary solid with the displacement of the
original solid. By Betti's theorem \cite{ghali2017structural}, the work done by the traction on auxiliary solid through the displacements on the
original solid is equal to the work done by the traction on original solid through the displacements on the auxiliary solid. Neither the traction on
original solid nor the displacements on the auxiliary solid depends on the geometry of the shape.
Thus, the first integral in Eq.\,\eqref{eq:E_suo} does not depend on the shape of the void, i.e. it is a constant $C$.  Thus,
\begin{equation}
    \Pi_{\mbox{el}} = E_{\mbox{int}} - W_{\mbox{ext}} = - E_{\mbox{int}} =\tfrac{1}{2} \int_{\partial D} n_i \sigma_{ij}^\infty
    u_j \, ds - C,
\end{equation}
which is Eq.\,\eqref{eq:12}.

%% file: Appendix_B.tex
\section{Derivation of Euler-Lagrange equation for void shape} \label{Appendix:B}
 
In this Appendix we derive the equilibrium condition for the surface by minimizing the total potential energy $\Pi$ given in Eq.~\eqref{eq:100} subject to the constraint of fixed void volume Eq.~\eqref{eq:area_constraint}.   

Expressing the geometry of the void in polar coordinates $r=r(\theta)$ and using the two-fold symmetry, the objective function for the total potential energy can be written as a functional of the shape $r(\theta)$ as
\begin{equation}
    \Pi = \Pi[r(\theta)]=\int^{\pi/2}_0\gamma[\omega(r,\dot{r})]\sqrt{r^2+\dot{r}^2}\,\mathrm{d}\theta
-\frac{1}{2}\int^{\pi/2}_0\int^{\infty}_r \epsilon_{ij}\sigma_{ij}R\,\mathrm{d}R\mathrm{d}\theta,
\end{equation}
where $\gamma[\omega(r,\dot{r})]$ represents the dependence of the anisotropic energy $\gamma(\omega)$ on the orientation angle $\omega$ via the shape function $r(\theta)$ through
\begin{equation}
\omega(r,\dot{r}) = \mbox{arccot}\left[ \frac{r \cos \theta + \dot{r} \sin \theta}{r \sin \theta - \dot{r} \cos \theta} \right] + \pi.
\end{equation}
Using a Lagrange multiplier $\mu$ to incorporate the constraint of constant void area into the functional we obtain a new functional
\begin{equation}
    H[r(\theta)]=\int_0^{\pi/2}L(r,\dot{r})\,\mathrm{d}\theta,
\end{equation}
where
\begin{equation}
    L(r,\dot{r})=\gamma[\omega(r,\dot{r})]\sqrt{{\dot{r}}^2+r^2}-
 \frac{1}{2}\int^{\infty}_r \epsilon_{ij}\sigma_{ij}R\,\mathrm{d}R
    -\mu\left(\frac{1}{2}r^2\,-\pi/4 \right),\label{eq:B4}
\end{equation}
The Euler-Lagrange equation for minimizing the functional is
\begin{equation}
    \frac{\partial L}{\partial r}- \frac{\mathrm{d}}{\mathrm{d}\theta}\left(\frac{ \partial L}{ \partial \dot{r}}\right)=0 \qquad \mbox{on} \quad 0<\theta<\pi/2.
\end{equation}
Denoting $\gamma' = d\gamma/d\omega$, the derivative terms in the Euler-Lagrange equation are
\begin{equation}
    \frac{\partial L}{\partial r}  =  \gamma'\cdot\frac{\dot{r}}{\sqrt{\dot{r}^2+r^2}}
    +\gamma\cdot\frac{r}{\sqrt{\dot{r}^2+r^2}}
    +\frac{1}{2}\epsilon_{ij}\sigma_{ij}r-\mu r,
\end{equation}
and
\begin{equation}
    \frac{\mathrm{d}}{\mathrm{d}\theta}\left(\frac{ \partial L}{ \partial \dot{r}}\right)= \frac{\gamma''(r^3+r^2\ddot{r}+2\dot{r}^2r)
        +\gamma'(\dot{r}^3+r^2\dot{r})+\gamma(\ddot{r}r^2-r\dot{r}^2)}
        {\left(\dot{r}^2+r^2\right)^{\frac{3}{2}}}.
\end{equation}
Thus the Euler-Lagrange equation for the surface can be reduced to
\begin{equation}
    (\gamma+\gamma'')\kappa 
    -\frac{1}{2}\epsilon_{ij}\sigma_{ij}+\mu =0 
    \qquad \mbox{on}\ L_1. \label{eq:112}
\end{equation}
where
\begin{equation}
\kappa = 
    \frac{(2\dot{r}^2-r\ddot{r}+r^2)}{(\dot{r}^2+r^2)^{3/2}}
    \label{eq:curvature}
\end{equation}
is the curvature of the surface ($\kappa>0$ for a convex void shape).

We now derive a simplified result for $\epsilon_{ij}\sigma_{ij}$ on $L_1$.
Inverting Hooke's law in a linear elastic solid Eq.\,\eqref{eq:constitutive_law} gives the relation between $\epsilon_{ij}$ and $\sigma_{ij}$ for the three-dimensional elasticity problem as
\begin{equation}
    \epsilon_{ij}=\frac{1}{E}\left[(1+\nu)\sigma_{ij}-\nu\delta_{ij}\sigma_{kk}\right], \label{eq:109}
\end{equation}
where $E$ and $\nu$ are Young's modulus and Poisson's ratio.  We consider plane strain elasticity for which $\epsilon_{33}=0$, then
\begin{equation}
    \sigma_{33}=\nu(\sigma_{11}+\sigma_{22}).
\end{equation}
The two-dimensional plane strain version of Eq.\,\eqref{eq:109} is thus
\begin{equation}
    \epsilon_{ij}=\frac{1+\nu}{E}\left[\sigma_{ij}-\nu\delta_{ij}(\sigma_{11}+\sigma_{22})\right]. \label{eq:111}
\end{equation}
Substitute for $\epsilon_{ij}$ using Eq.\,\eqref{eq:111} into $\epsilon_{ij}\sigma_{ij}$ to obtain
\begin{equation}
\epsilon_{ij}\sigma_{ij}  = \left[\sigma_x^2+\sigma_y^2-\nu(\sigma_x+\sigma_y)^2+2\tau_{xy}^2\right] 
=  \frac{1-\nu^2}{E}\left(\sigma_x+\sigma_y\right)^2 \quad \mbox{on}\ L_1,
\end{equation}
where we have utilized the stress-free condition on the inner boundary $L_1$ given by Eq.\,\eqref{eq:3}. Thus Eq.\,\eqref{eq:112} is equivalent to
\begin{equation}
(\gamma+\gamma'')\kappa-\frac{1-\nu^2}{2E}(\sigma_{xx}+\sigma_{yy})^2=-\mu
\quad \mbox{on}\ L_1.
\end{equation}
 
Finally, at the ends of the interval $\theta=0, \pi/2$, the natural boundary conditions associated with the Euler-Lagrange equations are that the boundary $r(\theta)$ is continuous and that the derivative of the shape $\dot{r}(\theta)$ must satisfy a jump conditions
\begin{gather}
    \left. \frac{\partial L}{ \partial \dot{r}} \right|^{\theta=0^+}_{\theta=0^-} =0, \\
    \left.\frac{\partial L}{ \partial \dot{r}}\right|^{\theta=\pi/2^+}_{\theta=\pi/2^-}=0.
\end{gather}
Using the assumed two-fold symmetry of the shape, the jump conditions reduce to
\begin{gather}
     \frac{\partial L}{ \partial \dot{r}}(\theta=0)=0, \\
    \frac{\partial L}{ \partial \dot{r}}(\theta=\pi/2)=0.
\end{gather}
Note that since the elasticity term in $L(r,\dot{r})$ given by Eq.\,\eqref{eq:B4} does not depend explicitly on $\dot{r}$, these boundary conditions are independent of the elastic energy and depend only on the surface energy contribution, which can be simplified to obtain the corner angle conditions 
\begin{equation}
\frac{\gamma'(\omega)}{\gamma(\omega)} = \frac{\dot{r}(\theta)}{r(\theta)}  \hs \mbox{at $\theta=0, \pi/2$}.
\label{corner_bc}
\end{equation}
Note also that for the case of isotropic surface energy, $\gamma'(\omega)$=0, the boundary conditions reduce to $\dot{r}=0$ at $\theta=0,\pi/2$, i.e.\ $r(\theta)$ has continuous first derivatives and no corner.  For the case of anisotropic surface energy $\gamma'(\omega)$ has nonzero values which can permit the existence of a corner with $\dot{r} \ne 0$ i.e.\ a discontinuous first derivative at the corner.

%% file: Appendix_C.tex
\section{Asymptotic analysis of corner singularity}\label{Appendix:C}

\setcounter{figure}{0}

\begin{figure}[!t]
    \centering
    \includegraphics[width=0.7\textwidth]{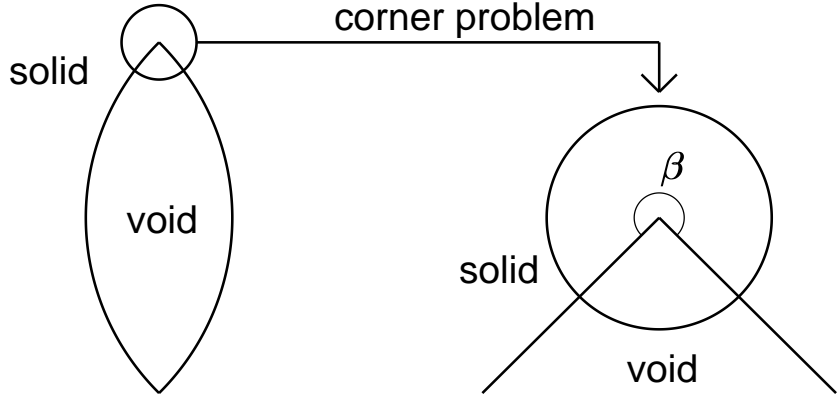}
    \caption{Boundary shape of the corner problem.}
    \label{fig:11}
\end{figure}

This section employs asymptotic analysis \cite{de1981asymptotic} to characterize the singularity in the elastic field near the corner of a void. The Goursat function $\varphi(z)$ is smooth along the boundary if the
shape of the void is smooth (with no corners). For a convex void
shape with a corner, $\varphi'(z)$ has a singularity at the corner and the
stress goes to infinity when $z$ approaches the corner \cite{savin1970stress,williams1952stress}.
We analyze the stress asymptotically to determine the order of the
singularity as a function of the corner angle. One would expect that in the 
vicinity of the corner the local problem
is that of a semi-infinite wedge (see Fig.\,\ref{fig:11}) with singular stresses as found in \cite{williams1952stress}. 

The asymptotic derivation of the wedge geometry local to the corner is as follows.
Setting the vertex of the corner as the origin, 
the biharmonic equation Eq.\,\eqref{eq:biharmonic} in polar coordinates is given in \cite{soutas2012elasticity01}
as
\begin{equation}
\left(\frac{\partial^{2}}{\partial r^{2}}+\frac{1}{r}\frac{\partial}{\partial r}+\frac{1}{r^{2}}\frac{\partial^{2}}{\partial\theta^{2}}\right)^{2}W=0,\label{eq:35}
\end{equation}
and the corresponding stress components in polar coordinates are given
by
\begin{gather}
\sigma_{rr}=\frac{1}{r}\frac{\partial W}{\partial r}+\frac{1}{r^{2}}\frac{\partial^{2}W}{\partial\theta^{2}}, \\
\sigma_{\theta\theta}=\frac{\partial^{2}W}{\partial r^{2}}, \\
\sigma_{r\theta}=-\frac{\partial}{\partial r}\left(\frac{1}{r}\frac{\partial W}{\partial\theta}\right),
\end{gather}
where $\sigma_{rr}$ is the stress in the radial direction, $\sigma_{\theta\theta}$
is the stress in the $\theta$ direction, and $\sigma_{r\theta}$ is the
shear stress.
Letting $r=\varepsilon\tilde{r}$ and $W(r,\theta)=\widetilde{W}(\tilde{r},\theta)$
with $\varepsilon\ll1$ to find the corner solution, the biharmonic
Eq.\,\eqref{eq:35} becomes
\begin{equation}
\frac{1}{\varepsilon^{2}}\left(\frac{\partial^{2}}{\partial\tilde{r}^{2}}+\frac{1}{\tilde{r}}\frac{\partial}{\partial\tilde{r}}+\frac{1}{\tilde{r}^{2}}\frac{\partial^{2}}{\partial\theta^{2}}\right)^{2}\widetilde{W}=0.\label{eq:36}
\end{equation}
Stresses of the corner problem are
\begin{gather}
\widetilde{\sigma_{\tilde{r}\tilde{r}}}=\frac{1}{\tilde{r}}\frac{\partial\widetilde{W}}{\partial\tilde{r}}+\frac{1}{\tilde{r}^{2}}\frac{\partial^{2}\widetilde{W}}{\partial\theta^{2}}=\varepsilon^{2}\sigma_{rr}, \\
\widetilde{\sigma_{\theta\theta}}=\frac{\partial^{2}\widetilde{W}}{\partial\tilde{r}^{2}}=\varepsilon^{2}\sigma_{\theta\theta}, \\
\widetilde{\sigma_{\tilde{r}\theta}}=-\frac{\partial}{\partial\tilde{r}}\left(\frac{1}{\tilde{r}}\frac{\partial\widetilde{W}}{\partial\theta}\right)=\varepsilon^{2}\sigma_{r\theta}.
\end{gather}
The differential equation of the corner problem is still a biharmonic
equation from Eq.\,\eqref{eq:36}.  Since the far-field conditions $\sigma_{rr},\sigma_{r\theta},\sigma_{\theta\theta}$
are finite away from the corner, we have $\widetilde{\sigma_{rr}},\widetilde{\sigma_{r\theta}},\widetilde{\sigma_{\theta\theta}}\rightarrow0$
as $\tilde{r}\rightarrow\infty$. 
Hence the local corner problem is identical to the wedge problem with no traction on the boundary.  The solution
to the wedge problem with such ``free-free'' boundary conditions is given by separation of variables in \cite{williams1952stress}. The
form of $\varphi$, $\varphi'$ and stresses near the corner are
\begin{equation}
 \varphi(r,\theta)\sim r^{\lambda-1}, \hspace{1em} \varphi'(r,\theta)\sim r^{\lambda-2}, \hspace{1em} \sigma_{rr},\sigma_{r\theta}, \sigma_{\theta\theta}\sim r^{\lambda-2},
\end{equation}
where $\lambda$ is the solution of  \cite{williams1952stress}
\begin{equation}
\sin\left[(\lambda-1)\beta\right]=-(\lambda-1)\sin\beta, \label{order_singularity}
\end{equation}
and where $\beta$ is the solid corner angle (see Fig.~\ref{fig:11}).
From Eq.\,\eqref{order_singularity}, we can determine the behavior of the stresses
near the corner. For the case of $\beta<\pi$, $\lambda$ is greater than
$2$ and there is no singularity near the corner. For the case of $\beta=\pi$,
the boundary is a straight line and the stresses are constant with $\lambda=2$.
When $\pi<\beta\leq 2\pi$ (see Fig.\,\ref{fig:11}), $\lambda$ is between
$3/2$ and $2$ which gives a stress singularity near the corner. As $\beta$
approaches $2\pi$, $\lambda$ approaches $3/2$. For all cases in which $\pi<\beta\le 2\pi$,
the singularity in the stresses near the corner is an integrable singularity.

Thus the stresses near the corner in our
problem have the same order singularity as the
wedge problem with the same corner angle. Note, however, that while 
the asymptotic analysis determines the order of the
weak singularity near the corner, it does not determine
the amplitude coefficient. Determination of the coefficient of the corner singularity
requires the full solution of the elasticity problem over the entire domain because
the weakly-singular corner term decays slowly with $r$ and thus contributes
to the elastic stress over the entire boundary of the void.

%% file: Appendix_D.tex
\section{Derivation of shape corner terms from corner terms of Goursat function} \label{Appendix:D}

Here we determine the asymptotic behavior of the void shape due to the influence of the elastic singularity near a corner. Since our elasticity formulation is in terms of the Goursat function $\varphi$, we determine how the shape $r(\theta)$ relates to $\varphi(\theta)$ near a corner.
We consider here the analysis near a corner at $\theta=0$. Similar results can be derived for a corner at $\theta=\pi/2$.
To find the relation between coefficients of $r$ and $\varphi$, we analyze the behavior of $\varphi'(z)$ near the corner at $\theta=0$. 
The shape of the surface near the corner is described by the equilibrium condition  Eq.\,\eqref{eq:114} given here as 
\begin{equation}
\left[\gamma(\omega)+\gamma''(\omega)\right]\kappa 
-\tfrac{1}{4}\Lambda
\left[ 1 + \chi + 4 \mbox{Re}\left\{\varphi'(z)\right\} \right]^2=\mu \quad \mbox{on}\ \partial D, \label{eq:D3}
\end{equation} 
Since Eq.\,\eqref{eq:D3} has a singular term $\varphi'(z)$, the only way the surface equilibrium equation can be satisfied is if there is a singularity
in the curvature $\kappa$ to cancel the elastic singularity.
From the expansion of $\varphi$ in our numerical method in 
Eq.\,\eqref{eq:115}
we have
\begin{equation}
\varphi'(\theta)\sim(\lambda_1-1)a_1\theta^{\lambda_1-2}+A+i\cdot((\lambda_1-1)b_1\theta^{\lambda_1-2}+B)
\qquad\mbox{as}\quad \theta\rightarrow 0 
\end{equation}
where $A$ and $B$ are known numerically-determined constants. 
Thus, balancing the elastic singularity of order $\lambda_1-2$ in 
\eqref{eq:D3} requires
a curvature term with singularities of order $\lambda_1-2$ and $2(\lambda_1-2)$.  From the definition of the curvature in terms 
of $r(\theta)$ from \eqref{eq:curvature} we deduce that
\begin{equation}
r''(\theta)\sim c_1(2\lambda_1-2)(2\lambda_1-3)\theta^{2\lambda_1-4}+c_2\lambda_1(\lambda_1-1)\theta^{\lambda_1-2}+\mathcal{O}(1)
\qquad\mbox{as}\quad \theta\rightarrow 0 \label{eq:B2}
\end{equation} 
where $c_1$ and $c_2$ are constants that must be chosen to satisfy \eqref{eq:D3}.

We evaluate the terms in Eq.\,\eqref{eq:D3} in the vicinity of $\theta=0$ to determine $c_1$ and $c_2$ of $r(\theta)$ in terms of the coefficients of
$\varphi(\theta)$. In polar coordinates,
\begin{equation}
\lim\limits_{\theta\rightarrow 0} z'(\theta)=r'(0)+ir(0),
\end{equation}
and by geometry
\begin{equation}
r'(0)/r(0)=-\tan({\alpha_1/2-\pi/2})=\cot(\alpha_1/2)
\end{equation}
where $\alpha_1$ is the interior corner angle. Then 
\begin{equation}
\begin{aligned}
\lim\limits_{\theta\rightarrow 0}\varphi'(z) & =\lim\limits_{\theta\rightarrow 0}\frac{\varphi'(\theta)}{z'(\theta)}
\sim\frac{(\lambda_1-1)a_1\theta^{\lambda_1-2}+A+i\cdot((\lambda_1-1)b_1\theta^{\lambda_1-2}+B)}{(\cot(\alpha_1/2)+i)r(0)}\\
& = \frac{(\lambda_1-1)(a_1+ib_1)(\cot(\alpha_1/2)-i)}{r(0)(\cot^2(\alpha_1/2)+1)}\theta^{\lambda_1-2}
+\frac{(A+iB)(\cot(\alpha_1/2)-i)}{r(0)(\cot^2(\alpha_1/2)+1)}.
\end{aligned}
\end{equation}
Taking the real part we obtain
\begin{equation}
\lim\limits_{\theta\rightarrow 0}\mbox{Re}\left\{\varphi'(z)\right\}\sim
\frac{(\lambda_1-1)(a_1\cot(\alpha_1/2)+b_1)\theta^{\lambda_1-2}+(A\cot(\alpha_1/2)+B)}{r(0)(\cot^2(\alpha_1/2)+1)}.
\end{equation}
Considering other terms in Eq.\,\eqref{eq:D3}, 
\begin{gather}
\lim\limits_{\theta\rightarrow 0}(\gamma+\gamma'')=1-15\varepsilon\cos(2\alpha_1), \\
\lim\limits_{\theta\rightarrow 0}\kappa\sim\lim\limits_{\theta\rightarrow 0}\frac{-r\cdot r''}{({r'}^2+r^2)^{3/2}}
=-\frac{c_1(2\lambda_1-2)(2\lambda_1-3)\theta^{2\lambda_1-4}+c_2\lambda_1(\lambda_1-1)\theta^{\lambda_1-2}}{{r(0)}^2(\cot^2(\alpha_1/2)+1)^{3/2}}.
\end{gather}
Comparing the coefficients of the $\theta^{2\lambda_1-4}$ term in Eq.\,\eqref{eq:D3} we obtain
\begin{equation}
c_1=-\frac{4\Lambda(\lambda_1-1)(a_1\cot(\alpha_1/2)+b_1)^2}{2(\cot^2(\alpha_1/2)+1)^{1/2}(1-15\varepsilon\cos(2\alpha_1))(2\lambda_1-3)}. 
\end{equation}
Comparing the coefficients of the $\theta^{\lambda_1-2}$ term in Eq.\,\eqref{eq:D3} we obtain
\begin{equation}
\begin{gathered}
c_2=-\frac{2\Lambda(1+\chi)(a_1\cot(\alpha_1/2)+b_1)(\cot^2(\alpha_1/2)+1)^{1/2}r(0)}{\lambda_1(1-15\varepsilon\cos(2\alpha_1))}\\
-\frac{8\Lambda(\lambda_1-1)(a_1\cot(\alpha_1/2)+b_1)(A\cot(\alpha_1/2)+B)}{2(\cot^2(\alpha_1/2)+1)^{1/2}(1-15\varepsilon\cos(2\alpha_1))\lambda_1}. 
\end{gathered}
\end{equation}
Following the same logic we can find the corner term coefficients at 
$\theta=\pi/2$.

%% file: DoesElasticStressModifyEquilbriumCornerAngle.bbl
\begin{thebibliography}{10}
\expandafter\ifx\csname url\endcsname\relax
  \def\url#1{\texttt{#1}}\fi
\expandafter\ifx\csname urlprefix\endcsname\relax\def\urlprefix{URL }\fi
\expandafter\ifx\csname href\endcsname\relax
  \def\href#1#2{#2} \def\path#1{#1}\fi

\bibitem{wenzel1936}
R.~N. Wenzel, Resistance of solid surfaces to wetting by water, Ind. Eng. Chem.
  28~(8) (1936) 988--994.
\newblock \href {https://doi.org/https://doi.org/10.1021/ie50320a024}
  {\path{doi:https://doi.org/10.1021/ie50320a024}}.

\bibitem{young1805}
T.~Young, An essay on the cohesion of fluids, Philos. Trans. R. Soc. 95 (1805)
  65--87.
\newblock \href {https://doi.org/https://doi.org/10.1098/rstl.1805.0005}
  {\path{doi:https://doi.org/10.1098/rstl.1805.0005}}.

\bibitem{choi2009}
W.~Choi, A.~Tuteja, J.~M. Mabry, R.~E. Cohen, G.~H. McKinley, A modified
  cassie--baxter relationship to explain contact angle hysteresis and
  anisotropy on non-wetting textured surfaces, J. Colloid Interface Sci.
  339~(1) (2009) 208--216.
\newblock \href {https://doi.org/https://doi.org/10.1016/j.jcis.2009.07.027}
  {\path{doi:https://doi.org/10.1016/j.jcis.2009.07.027}}.

\bibitem{srolovitz1989}
D.~J. Srolovitz, On the stability of surfaces of stressed solids, Acta Mater.
  37~(2) (1989) 621--625.
\newblock \href {https://doi.org/https://doi.org/10.1016/0001-6160(89)90246-0}
  {\path{doi:https://doi.org/10.1016/0001-6160(89)90246-0}}.

\bibitem{cabrera1964}
N.~Cabrera, The equilibrium of crystal surfaces, Surf Sci 2 (1964) 320--345.
\newblock \href {https://doi.org/https://doi.org/10.1016/0039-6028(64)90073-1}
  {\path{doi:https://doi.org/10.1016/0039-6028(64)90073-1}}.

\bibitem{Voorhees1984}
P.~Voorhees, S.~Coriell, G.~McFadden, R.~Sekerka, The effect of anisotropic
  crystal-melt surface tension on grain boundary groove morphology, J. Cryst.
  Growth 67~(3) (1984) 425--440.
\newblock \href {https://doi.org/https://doi.org/10.1016/0022-0248(84)90035-6}
  {\path{doi:https://doi.org/10.1016/0022-0248(84)90035-6}}.

\bibitem{wulff1901}
G.~Wulff, Zur frage der geschwindigkeit des wachsthums und der aufl{\"o}sung
  der krystallfl{\"a}chen, Z. Krystallog 34 (1901) 449--530.
\newblock \href {https://doi.org/https://doi.org/10.1524/zkri.1901.34.1.449}
  {\path{doi:https://doi.org/10.1524/zkri.1901.34.1.449}}.

\bibitem{spencer2004}
B.~J. Spencer, Asymptotic solutions for the equilibrium crystal shape with
  small corner energy regularization, Phys Rev E 69~(1) (2004) 011603.
\newblock \href {https://doi.org/https://doi.org/10.1103/PhysRevE.69.011603}
  {\path{doi:https://doi.org/10.1103/PhysRevE.69.011603}}.

\bibitem{williams1952stress}
M.~L. Williams, Stress singularities resulting from various boundary conditions
  in angular corners of plates in extension, J. Appl. Mech. 19 (1952) 526--528.
\newblock \href {https://doi.org/https://doi.org/10.1115/1.4010553}
  {\path{doi:https://doi.org/10.1115/1.4010553}}.

\bibitem{srolovitz2001}
D.~Srolovitz, S.~H. Davis, Do stresses modify wetting angles?, Acta Mater.
  49~(6) (2001) 1005--1007.
\newblock \href {https://doi.org/https://doi.org/10.1016/s1359-6454(01)00004-0}
  {\path{doi:https://doi.org/10.1016/s1359-6454(01)00004-0}}.

\bibitem{siegel2004}
M.~Siegel, M.~Miksis, P.~Voorhees, Evolution of material voids for highly
  anisotropic surface energy, J. Mech. Phys. Solids 52 (2004) 1319--1353.
\newblock \href {https://doi.org/https://doi.org/10.1016/j.jmps.2003.11.003}
  {\path{doi:https://doi.org/10.1016/j.jmps.2003.11.003}}.

\bibitem{oh1995method}
H.-S. Oh, I.~Babu, et~al., The method of auxiliary mapping for the finite
  element solutions of elasticity problems containing singularities, J. Comput.
  Phys. 121~(2) (1995) 193--212.
\newblock \href {https://doi.org/https://doi.org/10.1016/s0021-9991(95)90017-9}
  {\path{doi:https://doi.org/10.1016/s0021-9991(95)90017-9}}.

\bibitem{liu2012singular}
P.~Liu, T.~Bui, C.~Zhang, T.~Yu, G.~Liu, M.~Golub, The singular edge-based
  smoothed finite element method for stationary dynamic crack problems in 2d
  elastic solids, Comput Methods Appl Mech Eng 233 (2012) 68--80.
\newblock \href {https://doi.org/https://doi.org/10.1016/j.cma.2012.04.008}
  {\path{doi:https://doi.org/10.1016/j.cma.2012.04.008}}.

\bibitem{soutas2012potential}
R.~W. Soutas-Little, Elasticity, Dover Publications, New York, 2010, Ch.~11,
  pp. 251--263.

\bibitem{suo1994}
Z.~Suo, W.~Wang, Diffusive void bifurcation in stressed solid, J. Appl. Phys.
  76 (1994) 3410--3421.
\newblock \href {https://doi.org/https://doi.org/10.1063/1.357471}
  {\path{doi:https://doi.org/10.1063/1.357471}}.

\bibitem{bertotti1998hysteresis}
G.~Bertotti, Hysteresis in magnetism: For Physicists, Materials Scientists, and
  Engineers, Academic Press, 1998.
\newblock \href
  {https://doi.org/https://doi.org/10.1016/b978-0-12-093270-2.x5048-x}
  {\path{doi:https://doi.org/10.1016/b978-0-12-093270-2.x5048-x}}.

\bibitem{herring1951}
C.~Herring, Some theorems on the free energies of crystal surfaces, Phys. Rev.
  82~(1) (1951) 87.
\newblock \href {https://doi.org/https://doi.org/10.1103/physrev.82.87}
  {\path{doi:https://doi.org/10.1103/physrev.82.87}}.

\bibitem{pimpinelli1998physics}
A.~Pimpinelli, J.~Villain, Physics of crystal growth, Cambridge University
  Press, Cambridge, 1998, Ch.~3, pp. 50--55.
\newblock \href {https://doi.org/https://doi.org/10.1017/CBO9780511622526}
  {\path{doi:https://doi.org/10.1017/CBO9780511622526}}.

\bibitem{mikhlin1957integral}
S.~G. Mikhlin, Integral Equations, Pergamon Press, London, 1957.
\newblock \href {https://doi.org/https://doi.org/10.1016/C2013-0-08209-6}
  {\path{doi:https://doi.org/10.1016/C2013-0-08209-6}}.

\bibitem{wang2021}
W.~Wang, B.~J. Spencer, Numerical solution for the stress near a hole with
  corners in an infinite plate under biaxial loading, J. Eng. Math. 127~(13)
  (2021) 1--25.
\newblock \href {https://doi.org/https://doi.org/10.1007/s10665-021-10104-8}
  {\path{doi:https://doi.org/10.1007/s10665-021-10104-8}}.

\bibitem{gonzalez2008first}
O.~Gonzalez, A.~M. Stuart, A first course in continuum mechanics, Cambridge
  University Press, Cambridge, 2008.
\newblock \href {https://doi.org/https://doi.org/10.1017/CBO9780511619571}
  {\path{doi:https://doi.org/10.1017/CBO9780511619571}}.

\bibitem{Muskhelishvili}
N.~I. Muskhelishvili, Some basic problems of the mathematical theory of
  elasticity, Noordhoff, Groningen, 1953.
\newblock \href {https://doi.org/https://doi.org/10.1007/978-94-017-3034-1}
  {\path{doi:https://doi.org/10.1007/978-94-017-3034-1}}.

\bibitem{Hoskins2019}
J.~G. Hoskins, V.~Rokhlin, K.~Serkh, On the numerical solution of elliptic
  partial differential equations on polygonal domains, SIAM J. Sci. Comput.
  41~(4) (2019) A2552--A2578.
\newblock \href {https://doi.org/https://doi.org/10.1137/18m1199034}
  {\path{doi:https://doi.org/10.1137/18m1199034}}.

\bibitem{hesthaven2007spectral}
J.~S. Hesthaven, S.~Gottlieb, D.~Gottlieb, Spectral methods for time-dependent
  problems, Vol.~21, Cambridge University Press, Cambridge, 2007.
\newblock \href {https://doi.org/https://doi.org/10.1017/CBO9780511618352}
  {\path{doi:https://doi.org/10.1017/CBO9780511618352}}.

\bibitem{dennis1981adaptive}
J.~E. Dennis~Jr, D.~M. Gay, R.~E. Walsh, An adaptive nonlinear least-squares
  algorithm, ACM Trans. Math. Softw. 7~(3) (1981) 348--368.
\newblock \href {https://doi.org/https://doi.org/10.1145/355958.355966}
  {\path{doi:https://doi.org/10.1145/355958.355966}}.

\bibitem{MatlabOTB}
Matlab optimization toolbox, the MathWorks, Natick, MA, USA (Version 8.4).

\bibitem{ling1948stresses}
C.-B. Ling, The stresses in a plate containing an overlapped circular hole, J.
  Appl. Phys. 19 (1948) 405--411.
\newblock \href {https://doi.org/https://doi.org/10.1063/1.1715080}
  {\path{doi:https://doi.org/10.1063/1.1715080}}.

\bibitem{weir2008implications}
G.~Weir, Implications from the ratio of surface tension to bulk modulus and
  nearest neighbour distance, for planar surfaces, Proc. Math. Phys. Eng. Sci.
  464 (2008) 2281--2292.
\newblock \href {https://doi.org/https://doi.org/10.1098/rspa.2007.0360}
  {\path{doi:https://doi.org/10.1098/rspa.2007.0360}}.

\bibitem{ghali2017structural}
A.~Ghali, A.~M. Neville, T.~G. Brown, Structural analysis: a unified classical
  and matrix approach, CRC Press, Boca Raton, 2017.
\newblock \href {https://doi.org/https://doi.org/10.1201/b22004}
  {\path{doi:https://doi.org/10.1201/b22004}}.

\bibitem{de1981asymptotic}
N.~G. De~Bruijn, Asymptotic methods in analysis, Vol.~4, Dover Publications,
  New York, 1981.

\bibitem{savin1970stress}
G.~N. Savin, Stress distribution around holes, NASA, Washington, D.C., 1970.

\bibitem{soutas2012elasticity01}
R.~W. Soutas-Little, Elasticity, Dover Publications, New York, 2010, Ch.~8, pp.
  155--193.

\end{thebibliography}
